\documentclass[pdftex,12pt]{amsart}
\usepackage[dvips]{graphicx,color}
\usepackage[pdftex,colorlinks=true, pdfstartview=FitV, linkcolor=blue, citecolor=blue, urlcolor=blue]{hyperref}
\usepackage{a4wide}
\usepackage{amsthm}
\usepackage{amsopn}
\usepackage{amssymb}
\usepackage{amsmath}
\usepackage{latexsym}
\usepackage{graphics}
\usepackage[utf8]{inputenc}
\usepackage[T1]{fontenc}
\usepackage{amsfonts}
\usepackage{color}
\usepackage[all,cmtip]{xy}
\usepackage{mathrsfs}

\newtheorem{theorem}{Theorem}[section]
\newtheorem{lemma}[theorem]{Lemma}
\newtheorem{corollary}[theorem]{Corollary}
\newtheorem{proposition}[theorem]{Proposition}
\newtheorem{conjecture}{Conjecture}

\newenvironment{namedtheorem}[1]{\hfill\\ \noindent{\bf#1} \ \it}{\hfill\\}
\newtheorem{mainthm}{Theorem}

\theoremstyle{definition}

\newtheorem{definition}[theorem]{Definition}
\newtheorem{example}[theorem]{Example}

\newtheorem{remark}[theorem]{Remark}

\newcommand{\ble}{\begin{lemma}}
\newcommand{\ele}{\end{lemma}}
\newcommand{\bth}{\begin{theorem}}
\renewcommand{\eth}{\end{theorem}}
\newcommand{\bpr}{\begin{proposition}}
\newcommand{\epr}{\end{proposition}}
\newcommand{\bco}{\begin{corollary}}
\newcommand{\eco}{\end{corollary}}
\newcommand{\bcon}{\begin{conjecture}}
\newcommand{\econ}{\end{conjecture}}
\newcommand{\bde}{\begin{definition}}
\newcommand{\ede}{\end{definition}}
\newcommand{\bre}{\begin{remark}}
\newcommand{\ere}{\end{remark}}
\newcommand{\bex}{\begin{example}}
\newcommand{\eex}{\end{example}}
\newcommand{\beq}{\begin{equation}}
\newcommand{\eeq}{\end{equation}}

\newcommand{\pf}{\noindent{\bf Proof}\hspace{7pt}}

\newcommand{\bpf}{\pf}
\newcommand{\epf}{\qed}
\newcommand{\into}{\hookrightarrow}

\newcommand{\sbe}{\subseteq}
\newcommand{\sbne}{\subsetneq}

\newcommand{\spe}{\supseteq}
\newcommand{\De}{\Delta}

\newcommand{\LatU}{\Upsilon}
\newcommand{\Lat}{Y}
\newcommand{\cA}{{\mathcal A}}
\newcommand{\cC}{{\mathcal C}}
\newcommand{\cD}{{\mathcal D}}
\newcommand{\cE}{{\mathcal E}}

\newcommand{\cL}{{\mathcal L}}
\newcommand{\cO}{{\mathcal O}}
\newcommand{\cP}{{\mathscr P}}
\newcommand{\A}{{\mathsf{A}}}

\newcommand{\amG}{{\mathscr{G}}}
\newcommand{\amB}{{\mathscr{B}}}
\newcommand{\geomI}{{\mathscr{I}}}
\newcommand{\Dt}{{\mathbf{D}}}
\newcommand{\GUD}{{\mathbf{GUD}}}
\newcommand{\SUD}{{\mathbf{SUD}}}
\newcommand{\GD}{{\mathbf{GD}}}

\newcommand{\F}{{\mathsf{F}}}
\newcommand{\G}{{\mathbf{G}}}
\newcommand{\Q}{{\mathsf{Q}}}
\newcommand{\hG}{{\mathbf{\hat G}}}

\newcommand{\Inv}{\mathop{\rm Inv}\nolimits}
\renewcommand{\L}{{\mathbf{L}}}
\newcommand{\amL}{{\mathscr{L}}}
\newcommand{\n}{{\mathsf{n}}}
\newcommand{\M}{M}
\newcommand{\R}{{\mathsf{R}}}
\newcommand{\U}{U}
\newcommand{\V}{V}

\newcommand{\kaut}{\alpha}
\DeclareMathOperator{\Z}{Z}
\newcommand{\tG}{\tilde{\G}}
\newcommand{\tGamma}{\tilde{\Gamma}}
\newcommand{\tphi}{\tilde{\phi}}

\newcommand{\Cyc}{C}

\newcommand{\adj}{\diamond}
\DeclareMathOperator{\Aut}{Aut}
\DeclareMathOperator{\diag}{diag}
\DeclareMathOperator{\Det}{Det}

\DeclareMathOperator{\End}{End}
\DeclareMathOperator{\id}{id}
\DeclareMathOperator{\PSL}{PSL}
\DeclareMathOperator{\GU}{GU}
\DeclareMathOperator{\SU}{SU}
\DeclareMathOperator{\SL}{SL}
\DeclareMathOperator{\GL}{GL}
\DeclareMathOperator{\PG}{PG}
\newcommand{\FF}{{\mathbb F}}
\newcommand{\NN}{{\mathbb N}}
\newcommand{\PP}{{\mathbb P}}
\newcommand{\ZZ}{{\mathbb Z}}
\DeclareMathOperator{\typ}{typ}
\newcommand{\Gm}{\Gamma}
\DeclareMathOperator{\opp}{opp}
\newcommand{\fk}{{\mathsf{k}}}
\newcommand{\dfn}{\em}
\newcommand{\after}{\mathbin{ \circ~}}

\DeclareMathOperator{\Sp}{Sp}

\DeclareMathOperator{\PGL}{PGL}
\DeclareMathOperator{\proj}{proj}

\newcommand{\sP}{{\mathsf{P}}}

\newcommand{\sA}{{\mathsf{A}}}
\newcommand{\tA}{\widetilde{A}}
\newcommand{\tI}{\widetilde{I}}
\newcommand{\tJ}{{\widetilde{J}}}
\DeclareMathOperator{\aff}{aff}
\DeclareMathOperator{\im}{im}

\newcommand{\vep}{\varepsilon}
\newcommand{\valuation}{v}
\newcommand{\mn}{\ \medskip \newline }
\newcommand{\ead}{\email}
\newenvironment{keyword}{Keywords: \keywords}{}

\newcommand{\MSC}[1]{MSC #1:\subjclass}
\newcommand{\sep}{\hspace{1ex}}

\begin{document}

\title
{Curtis-Tits groups generalizing Kac-Moody groups of type $\tA_{n-1}$}

\author{Rieuwert J. Blok}
\ead{blokr@member.ams.org}
\address{Department of Mathematics and Statistics\\
Bowling Green State University\\
Bowling Green, oh 43403\\
U.S.A.}
\author{Corneliu G. Hoffman}
\ead{C.G.Hoffman@bham.ac.uk}
\address{University of Birmingham\\
Edgbaston, B15 2TT\\
U.K.}

%
\begin{abstract}
In~\cite{BloHof2013} we define a Curtis-Tits group as a certain generalization of a Kac-Moody group.
We distinguish between orientable and non-orientable Curtis-Tits groups and identify all orientable Curtis-Tits groups as
 Kac-Moody groups associated to twin-buildings. 
 
In the present paper we construct all orientable as well as non-orientable Curtis-Tits groups with diagram $\tA_{n-1}$ ($n\ge 4$) over a field $\fk$ of size at least $4$. 
The resulting groups are quite interesting in their own right. The
orientable ones are related to Drinfeld's construction of vector
bundles over a non-commutative projective line and to the  classical groups over cyclic algebras. The non-orientable ones are related to expander graphs~\cite{BloHofVdo2012a} and have symplectic, orthogonal and unitary groups as quotients. 
\end{abstract}
\maketitle

%
\section{Introduction}
The theory of the infinite dimensional Lie algebras called Kac-Moody algebras was initially developed by Victor Kac and Robert Moody. The development of a  theory of Kac-Moody groups as analogues of Chevalley groups was made possible by the work of Kac and   Peterson. 
In~\cite{Ti1992} J. Tits gives an alternative definition of  a group of Kac-Moody type  as being a group with a twin-root datum, which implies that  they are symmetry groups of Moufang twin-buildings.

In~\cite{AbMu1997} P.~Abramenko and B.~M\"uhlherr generalize a celebrated theorem of Curtis and Tits on groups with finite BN-pair~\cite{Cur1965a,Ti1974} to groups of Kac-Moody type.
This theorem states that a Kac-Moody group $\G$ is the universal completion of an amalgam of rank two (Levi) subgroups, as they are arranged inside $\G$ itself. This result was later refined by Caprace~\cite{Cap2007}.
Similar results on Curtis-Tits-Phan type amalgams have been obtained in~\cite{BenGraHof2003,BenGraHof2007,BeSh2004,BlHo2008,BlHo2009,GraHofNic2005,GraHofShp2003,GraHorNic2006,GraHorNic2007}.
For an overview of that subject see K\"ohl~\cite{Gra2009}.

In order to describe the main result from~\cite{BloHof2013} we introduce some notation.
Let $\fk$ be a (commutative) field of order at least $4$.
Let $\Gamma$ be a connected simply-laced Dynkin diagram over an index set $I$ without triangles.
For any $J\sbe I$, let $\Gamma_J$ be the subdiagram supported by the node set $J$.
In~\cite{BloHof2013} we take the Curtis-Tits type results as a starting point and define a {\em Curtis-Tits amalgam} with diagram $\Gamma$ over $\fk$ to be an amalgam of groups such that the sub-amalgam corresponding to a two-element subset $J\sbe I$ is the amalgam of derived groups of standard Levi subgroups of some rank-$2$ group of Lie type $\Gamma_J$ over $\fk$.
There is no a priori reference to an ambient group, nor to the existence of an associated (twin-) building. Indeed, there is no a priori guarantee that the amalgam will not collapse. Also, this definition clearly generalizes to other Dynkin diagrams.

We then classify all Curtis-Tits amalgams with diagram $\Gamma$ over $\fk$ using the following data (for similar results in special cases see~\cite{Dun2005,Gra2004}).
Viewing $\Gamma$ as a graph, for $i_0\in I$, let $\pi(\Gamma,i_0)$ denote the (first) fundamental group of  $\Gamma$ with base point $i_0$.
Also we let the group $\Aut(\fk)\times\langle \tau\rangle$ (with $\tau$ of order $2$) act as a subgroup of the stabilizer in $\Aut(\SL_2(\fk))$ of a fixed torus in $\SL_2(\fk)$;  $\tau$ denotes the transpose-inverse map with respect to that torus.
The main result of~\cite{BloHof2013} is the following.

\begin{namedtheorem}{Classification Theorem}\label{thm:CT amalgams} 
There is a natural bijection between isomorphism classes of Curtis-Tits amalgams with diagram $\Gamma$ over the field $\fk$ and group homomorphisms $\Theta\colon\pi(\Gamma, i_0)\to \langle \tau\rangle\times\Aut(\fk)$.
 \end{namedtheorem}

We call amalgams corresponding to homomorphisms $\Theta$ whose image lies inside $\Aut(\fk)$ ``orientable''; others are called ``non-orientable''.
It is not at all immediate that all non-orientable amalgams arising from the Classification Theorem are non-collapsing, i.e.\ that their universal completion is non-trivial.
We shall call a non-trivial group a {\em Curtis-Tits group} if it is the universal completion of a Curtis-Tits amalgam.
It is shown that orientable Curtis-Tits amalgams are precisely those arising from the Curtis-Tits theorem applied to a group of Kac-Moody type. Thus, groups of Kac-Moody type are orientable Curtis-Tits groups.

\subsection{Main results}\label{subsection:main results}
We now specify $\Gamma$ to be the Dynkin diagram of type $\tA_{n-1}$ labeled cyclically with index set $I=\{1,2,\ldots,n\}$, where $n\ge 4$.
The purpose of the present paper is to construct all orientable and non-orientable Curtis-Tits groups over $\fk$ with diagram $\Gamma$ and to study their properties. 

The paper is structured as follows.  In Section~\ref{section:ct-groups} we introduce the relevant notions
about amalgams and  describe all possible Curtis-Tits amalgams of type $\Gamma$ over $\fk$.
For each $\delta\in \Aut(\fk)\times\langle\tau\rangle$ we introduce a Curtis-Tits amalgam $\amG^\delta$ corresponding to $\delta$ via $\Theta$ as in the Classification Theorem and denote its universal completion $(\tG^\delta,\tphi^\delta)$.
In Section~\ref{section:orientable} we exhibit a non-trivial completion for orientable Curtis-Tits groups using a description of the corresponding twin-building.
In order to state the main result of this section we introduce the following notation.
For $\kaut\in \Aut(\fk)$, let $\R_\kaut =\fk\{t, t^{-1}\}$ be the ring of skew Laurent polynomials with coefficients in the field $\fk$ such that for $x\in  \fk$ we have  $txt^{-1}=x^\kaut$.
Let $\fk_\kaut$ be the fixed field of $\kaut$ in $\fk$.
 We use the Dieudonn\'e determinant to identify $\SL_n(\R_\kaut)$.
As usual, the center of a group $X$, is denoted $\Z(X)$.
We obtain the following.

\begin{mainthm}\label{mainthm:orientable CT group}
For $\kaut\in \Aut(\fk)$, the universal completion $\tG^\kaut$ of  $\amG^\kaut$ is an extension of $\SL_n(\R)$ by a subgroup $H$ of the center 
$Z(\tG^\kaut)$, which is isomorphic to a subgroup of 
 $\fk_\kaut^*$.

\end{mainthm}

In Section~\ref{section:non-orientable} we consider the case $\delta=\kaut\tau$ for some $\kaut\in \Aut(\fk)$ and exhibit a non-trivial completion of $\amG^\delta$.
Via Proposition~\ref{prop:non-oriented CT amalgam}  we obtain the first two parts of Theorem~\ref{mainthm:omega kac-moody construction} below. Demonstrating the universality and identification of the completion is more involved this time and takes up Subsections~\ref{subsection:flag transitive},~\ref{sec:simple-connectedness},~\ref{subsec:filtration}~and~\ref{subsec:proof of theorem 3}.

In order to state the main result of Section~\ref{section:non-orientable}, we introduce the following notation.
Let  
 $\sigma$ be the automorphism of $\R_{\alpha^2}$ inducing $\kaut^{-1}$ on $\fk$  and interchanging $t$ and $t^{-1}$ and let $\beta$ be the asymmetric $\sigma$-sesquilinear form on the free $\R_{\alpha^2}$-module $M$ with ordered basis $(e_1,\ldots,e_n,f_1,\ldots,f_n)$ having $\beta$-Gram matrix
 \begin{align}
B=\left(  \begin{array}{c|c} 0_n &  I_n \\ \hline  tI_n & 0_n \end{array}\right)\in \GL_{2n}(\R_{\alpha^2}).
\end{align}
 
 \begin{mainthm}\label{mainthm:omega kac-moody construction}
The group $\SU_{2n}(\R_{\alpha^2})$ of symmetries in $\SL_{2n}(\R_{\alpha^2})$ of the $\sigma$-sesquilinear form $\beta$ contains a completion of $\amG^\delta$.

Now suppose, in addition, that $|\fk|\ge 7$, that $\kaut\tau$ has finite order $s$, that $\fk/\fk_\kaut$ is a cyclic Galois extension and that the norm $N_{\fk_{\kaut^2}/\fk_\kaut}$ is surjective.
Then, the universal completion $\tG^\delta$ of  $\amG^\delta$ is an extension of $\SU_{2n}(\R_{\kaut^2})$ by a subgroup $H$ of the center 
$Z(\tG^\delta)$, which is isomorphic to a subgroup of the kernel of $N_{\fk_{\kaut^2}/\fk_\kaut}$.
\end{mainthm}
Finally, we note that some of these groups have been studied in a different context, namely that of abstract involutions of Kac-Moody groups~\cite{GraHorMuh2011}.  There, connectedness, but not simple-connectedness, of geometries such as those defined in Section~\ref{section:non-orientable} is proved.


\subsection{Applications: the orientable Curtis-Tits groups \texorpdfstring{$\SL_n(\R_{\kaut})$}{}}\label{sec:orient-curt-tits}


Let $\De=((\De_+,\delta_+),(\De_-,\delta_-), \delta_*)$ be the twin-building associated to the 
Kac-Moody group $\SL_n(\R_{\kaut})$. 
Then, the pairs of maximal residues from $\De_+$ and $\De_-$ that are opposite for the twinning correspond to vector bundles over the non-commutative projective line $\PP^1(\kaut)$ in the sense of Drinfel'd. 
More precisely, let $\fk\{t\}, \fk\{t^{-1}\}\le \R_{\kaut}$ be the
corresponding skew polynomial rings and fix $\mathbf{M}$ a free
$\R_{\kaut}$ module of rank $n$. Following
\cite{LauRapStu1993} and  \cite{Tae2007} one can define a rank $n$ vector bundle
over the non-commutative projective line $\mathbb{P}^1(\kaut)$ as a
collection $(M_+,M_-, \phi_+, \phi_-)$ where $M_\vep$ is a free $n$-dimensional module over $\fk\{t^{\vep 1}\}$ and
$\phi_\vep\colon M_\vep\otimes\R_\kaut \to \mathbf{M}$ is an isomorphism of $\R_\kaut$-modules.
By analogy to the commutative case (see  \cite{Ro1992a, Ro1992} for
example) one could describe the building structure in terms of these
vector bundles. We intend to explore these relations to number theory
in a future paper.

To give a different perspective on these groups we note that the skew Laurent polynomials are  closely related to cyclic algebras as defined by
Dickson. More precisely let $\fk'\le \fk$ be  a cyclic field extension, of degree $n$, and let $\kaut$ be the generator of its Galois group. Given any $b\in \fk'$, define  the $\fk'$-algebra $(\fk/\fk', \kaut, b)$ to be generated by the elements of $\fk$, viewed as an extension of $\fk'$, together with some
element $u$ subject to the following relations:
$$u^n = b, xu = ux^\kaut
 \mbox{ for } x \in \fk.$$

These algebras are central simple algebras.
A theorem due to Albert, Brauer, Hasse and Noether
~\cite{AlbHas1932,BrHaNo1932} says  that every central division algebra over a number field $\fk'$ is isomorphic to $(\fk/\fk', \kaut, b)$ for some $\fk, b, \kaut$. 
One constructs the map  $\epsilon_b:\R_{\kaut}\to (\fk/\fk', \kaut, b)$ via $t^{-1}\mapsto u$. This induces
a map $\epsilon_b\colon\SL_n(\R_{\kaut})\to \SL_n((\fk/\fk', \kaut, b))$, realizing
the linear groups over cyclic algebras as completions of the
Curtis-Tits amalgams.

\subsection{Applications: the purely non-orientable groups \texorpdfstring{$\G^\tau$}{}}
\label{sec:non-orient-groups}
We consider the situation described in Theorem~\ref{mainthm:omega kac-moody construction}, where we set $\delta=\tau$ (that is, $\kaut=\id_\fk$).
Then, $\R=\R_{\id}=\fk[t,t^{-1}]$ is the ring of  Laurent polynomials in the commuting variable $t$ over the field $\fk$.

It turns out that the group $\G^\tau=\SU_{2n}(\R,\beta)$ has some very interesting
natural quotients.
Let $\bar{\fk}$ denote the algebraic closure of $\fk$.
For any $b\in \bar{\fk}^*$ consider the specialization map
$\epsilon_b\colon\fk[t,t^{-1}]\rightarrow \bar{\fk}$ given by
$\epsilon_b(f)=f(b)$. The map induces a homomorphism $\epsilon_b\colon\SL_{2n}(\R) \to \SL_{2n}(\fk(b)) $. In some instances the map $b\leftrightarrow b^{-1}$ defines an automorphism of $\fk(b)$ and so one can define a
map $\epsilon_b\colon \G^\tau \to \SL_{2n}(\bar{\fk}) $

 The most important specialization maps are those given by evaluating $t$ at $b=\pm1$ or $b=\zeta$, a $(q^m+1)$-st root of $1$ where $q$ is a power of the characteristic. 

Consider first $b=-1$. In this case the automorphism $\sigma$ is trivial. 
Note that   for $g \in \G^\tau$ we have $\epsilon_{-1}(g)\in
\Sp_{2n}(\fk)$. In this case, the image of the group $\G^\tau$ is the group generated by the Curtis-Tits amalgam $\amL^\tau$ inside $\Sp_{2n}(\fk)$. 

Similarly, if $b=1$, the automorphism $\sigma$ is trivial and  the map $\epsilon_1$ takes $\G^\tau$ into $\Omega^+_{2n}(\fk)$.

Finally assume that  $\fk=\FF_q$ and $b \in \bar{\FF}_q$ is a primitive $(q+1)$-st root of $1$.
The $\FF_q$-linear map $\FF_q(b)\to \FF_q(b)$ induced by $\sigma$ sends $b$ to $b^{-1}$. 
Thus, $\sigma$ coincides with 
the Frobenius  automorphism of the field
$\FF_q(b)=\FF_{q^2}$.  
It is easy to verify that a change of
coordinates $e'_i=e_i$ and $f'_i =b f_i$ where $c^2 =b$ 
standardizes the Gram matrix of $\beta\after(\epsilon_b\times\epsilon_b)$ to a hermitian one, thus identifying the
image of $\epsilon_b$ with a subgroup of a conjugate of the unitary group $\SU_{2n}(q)$.
In~\cite{BloHofVdo2012a} it is shown that the image of this map is isomorphic to $\SU_{2n}(q)$.  This easily generalizes to the case where $b$ is a $(q^m+1)$-st root of unity and indeed to other cases where $a$ is Galois-conjugate to $b^{-1}$.
Also in~\cite{BloHofVdo2012a} we have shown that Cayley graphs of these groups form families of expander graphs.


\medskip
\paragraph{\bf Acknowledgement}
{This project was started during a visit to the Banff International Research Station and an earlier version was finished during a visit to the Mathematisches Forschungsinstitut Oberwolfach in the Research in Pairs program from October 25 until November 7, 2009. We thank both institutes for providing such a pleasant and stimulating research environment.}
We would also like to thank the anonymous referee for his/her careful reading of the manuscript and numerous suggestions for improvements. In particular, the proofs of Lemmas~\ref{lem:Wtau}~and~\ref{lem:sws'<>w} were significantly shortened as a result.

\section{Curtis-Tits groups}\label{section:ct-groups}
In this section we briefly recall the notion of a Curtis-Tits amalgam with diagram $\Gamma$ over $\fk$ from~\cite{BloHof2013}.
Recall that $\Gamma$ is the Dynkin diagram of type $\tA_{n-1}$ with nodes labeled cyclically by the elements of the index set $I=\{1,2,\ldots,n\}$ and that $\fk$ is a commutative field of order at least $4$.
\medskip

\bde\label{dfn:amalgam}
An {\em amalgam} over a poset $(\cP,\prec)$ is a collection $\amG=\{\G_x\mid x\in \cP\}$ of groups, together with a collection $\varphi=\{\varphi_x^y\mid x\prec y, x,y\in \cP\}$ of monomorphisms $\varphi_x^y\colon \G_x\into \G_y$, called {\em inclusion maps} such that whenever $x\prec y\prec z$, we have
 $\varphi_x^z=\varphi_y^z\after\varphi_x^y$.
A {\em completion} of $\amG$ is a group $\G$ together with a collection  $\phi=\{\phi_x\mid x\in \cP\}$ of homomorphisms $\phi_x\colon \G_x\to \G$, whose images generate $\G$, such that for any $x, y\in \cP$ with $x\prec y$ we have 
$\phi_y\after\varphi_x^y=\phi_x$.
The amalgam $\amG$ is {\em non-collapsing} if it has a non-trivial completion.
A completion $(\widetilde{\G},\tilde{\phi})$ is called {\em universal} if for any completion $(\G,\phi)$ there is a unique surjective group homomorphism $\pi\colon \widetilde{\G}\to \G$ such that $\phi=\pi\after\tilde{\phi}$. 
\ede

Before we define the Curtis-Tits amalgam $\amG^\delta$ we specify an action of the group $\Aut(\fk)\times\langle \tau\rangle$ (with $\tau$ of order $2$) on $\SL_2(\fk)$.
We let $\alpha\in \Aut(\fk)$ act entry-wise on $A\in \SL_2(\fk)$
and let $\tau$ act by sending each $A\in \SL_2(\fk)$ to its transpose inverse ${}^tA^{-1}$ with respect to the standard basis. Note that $\tau$ acts as an inner automorphism. 

\paragraph{Indexing convention}
Throughout the paper we shall adopt the following indexing conventions. Indices from $I$ shall be taken modulo $n$. For any $i\in I$, we set $(i)=I-\{i\}$.
Also subsets of $I$ of cardinality $1$ or $2$ appearing in subscripts are written without set-brackets.

\bde\label{dfn:standard CT structure}\label{dfn:standard pair}
Let $\cP=\{J\mid \emptyset\ne J\sbe I \mbox{ with }|J|\le 2\}$ and 
 $\prec$ denoting inclusion.
Given an element $\delta\in \Aut(\fk)\times\langle \tau\rangle$
the {standard universal} Curtis-Tits amalgam with diagram $\Gamma$
over $\fk$ corresponding to $\delta$ is the amalgam 
$\amG^\delta=\{\G_{i},\G_{i, j}, \psi_{i,j} \mid  i, j \in I\}$ 
over $\cP$, where, for every $i,j\in I$, we write $\psi_{i,j}=\psi_{\{i\}}^{\{i,j\}}$. Note that, due to our subscript conventions, we write  $\G_{i}=\G_{\{i\}}$  and $\G_{i,j}=\G_{\{i,j\}}$, 
where 
\begin{enumerate}
  \item[(SCT1)]  for any vertex $i$, we set $ \G_{i} = \SL_2(\fk)$ and for each pair $i,j \in I$,
 \begin{align*}\G_{i,j}\cong\begin{cases}
 \SL_3(\fk) & \mbox{if} \ \{i,j\}=\{i, i+1\} \\ \G_{i}{\times}
 \G_{j} & \mbox{if} \ \{i,j\}\ \ne \{i,i+1\} \end{cases},\end{align*}
  \item[(SCT2)] for $i=1,2,\ldots,n-1$ we have 
  $$ 
  \begin{array}{cc}
  \begin{array}{rl}
\psi_{i,i+1}\colon \G_{i}&\to \G_{i,i+1}\\
A & \mapsto  \left(\begin{array}{@{}cc@{}}A & 0 \\0 & 1\end{array}\right)
\end{array}
&
 \begin{array}{rl}
\psi_{i+1,i}\colon \G_{i+1}&\to \G_{i,i+1}\\
A & \mapsto  \left(\begin{array}{@{}cc@{}}1 & 0 \\0 & A\end{array}\right)
\end{array}
\end{array},$$
and we have 
  $$ 
  \begin{array}{cc}
  \begin{array}{rl}
\psi_{n,1}\colon \G_{n}&\to \G_{n,1}\\
A & \mapsto  \left(\begin{array}{@{}cc@{}}A & 0 \\0 & 1\end{array}\right)
\end{array}
&
 \begin{array}{rl}
\psi_{1,n}\colon \G_{1}&\to \G_{1,n}\\
A & \mapsto  \left(\begin{array}{@{}cc@{}}1 & 0 \\0 & A^{\delta}\end{array}\right)
\end{array}
\end{array},$$
whereas for all other pairs $(i,j)$, $\psi_{i,j}$ is the natural inclusion of $\G_{i}$ in $\G_{i}{\times}
\G_{j}$.
\end{enumerate}
We shall adopt the following shorthand:
$\G^+_{i}=\psi_{i,i+1}(\G_{i})$, 
$\G^-_{i}=\psi_{i,i-1}(\G_{i})$, 
where indices are taken modulo $n$.
\ede
By~\cite{BloHof2013}, every {universal} Curtis-Tits amalgam with Dynkin diagram $\tA_{n-1}$ over $\fk$
is isomorphic to $\amG^\delta$ for a unique $\delta\in\Aut(\fk)\times\langle\tau\rangle$.
We have chosen our setup such that $\amG^{\id}$ is the amalgam resulting from applying the Curtis-Tits theorem to the split Kac-Moody group $\SL_n(\fk[T,T^{-1}])$ of type $\tA_{n-1}$ with respect to its standard twin $BN$-pair. 

Note that the CT-amalgam $\amG^{\delta}$ has property $(D)$ as in \cite{BloHof2013}, that is, for any $i$ there exists a torus $D_{i} \in G_{i}$ so that 
\begin{align*}
\psi_{i,i+1}(D_{i})=N_{\G_{i}^+}(\G_{i+1}^-), & \mbox{ and }\\
\psi_{i,i-1}(D_{i})=N_{\G_{i}^-}(\G_{i-1}^+).
\end{align*}

\bde\label{dfn:OCT}
Note that since $|\fk|\ge 4$, a maximal split torus in $\SL_{2}(\fk)$ uniquely determines a pair of opposite root groups $X_{+}$ and $X_{-}$. We now choose one root group  $X_{i}$ normalized by the torus $D_{i}$ of $\G_{i}$ for each $i$. 
An {\it orientable Curtis-Tits (OCT) amalgam} (respectively orientable Curtis-Tits (OCT) group) is a non-collapsing Curtis-Tits amalgam that admits a system $\{X_{i}\mid i\in I\}$ of root groups as above such that for any $i,j\in I$, the groups $\psi_{i,j} (X_{i})$ and 
$\psi_{j,i}(X_{j})$ are contained in a common Borel subgroup  of $\G_{i,j}$. 
By the classification result in~\cite{BloHof2013} the amalgam $\amG^{\delta}$ is orientable if and only if $\delta \in \Aut(k)$.
\ede

In the remainder of this section we fix $\delta$ and we drop the superscript $\delta$, if no confusion arises.

Our methods are building theoretic and, for that reason we will need a thick version of a CT amalgam. To that end we need some notations. For any non-empty $J\sbe  I$ define the amalgam 
\begin{align*}
\amG_{J}=\{\G_{k}, \G_{k,l}, \psi_{k,l} | k,l \in J, k\ne l \}
\end{align*}
and let $(\G_J,\phi_J)$ be its universal completion.
Note that for $|J|\le 2$, $\G_J$ is the group from $\amG$ itself.
\ble\label{lem:cA_k}\label{lem:G_J}
Let $J\sbne I$ and let $J=\cup_{i}J_{i}$ be a decomposition of $J$ corresponding to connected components of the diagram $\Gm_J$ induced on the node set $J$.
Then $\G_{J}\cong 
\oplus_{i} \SL_{n_{i}+1}(\fk)$ where $|J_{i}|=n_i$. 
\ele
\bpf
For each $i$, we see that $\amG_{J_i}$ is exactly the unique Curtis-Tits amalgam of $\SL_{n_i+1}(\fk)$.
The result now follows from the Curtis-Tits theorem~\cite[Theorem 1]{Tim1998} (see also~\cite{Cur1965a,Ti1974,Tim2003,Tim2004,Tim2006}) recalling that $\SL_{n_i+1}(\fk)$ is the universal Chevalley group of type $A_{n_i}$ over $\fk$.
\epf
For any $m<n$, define an amalgam
\begin{align*}
\amG_{|m|} = \{\G_{J}, \psi_J^K\mid \emptyset\ne J\sbe K\sbne I, |K|\le m \},
\end{align*} 
where $\psi_J^K$ is given by universality.
We have 
 $\amG =\amG_{|2|}\sbe \amG_{|n-1|}$.
 
Recall that $(\tG,\tphi)$ is the universal completion of $\amG$. 
Let $(\tG_{|n-1|},\tphi_{|n-1|})$ be the universal completion of $\amG_{|n-1|}$. By construction of $\tG_{|n-1|}$, we have a  non-trivial map $\amG \to \tG_{|n-1|}$, so $\tG_{|n-1|}$ is a completion of $\amG$ and we get a surjective map $f\colon \tG\twoheadrightarrow \tG_{|n-1|} $.  Conversely, let $\emptyset\ne J\sbne I$. Then, the group $\tG_{J}=\langle \G_{i}| i \in J \rangle_{\tG}$ is a completion of the amalgam $\amG_{J}$ in $\tG$ and so there is a map $\G_{J} \mapsto \tG_{J}$. This means that $\tG$ is a completion of the amalgam $\amG_{|n-1|}$  and so there is a surjective map $g\colon \tG_{|n-1|}\to\tG$.
One now verifies that $g\after f\after\tphi_i=\tphi_i$ for all $i\in I$. By universality $g\after f$ is the identity map on $\tG$.
We have proved that
\ble\label{lem: fat amalgam}
$\amG$ and $\amG_{|n-1|}$ have the same universal completions.
\ele

We need to enlarge the amalgam even more. Consider $\G^\delta$ a completion of $\amG_{|n-1|}$. Denote by $\L_J$, respectively $\Dt_{i}$ the image of $\G_J$ respectively $D_{i}$ in $\G^\delta$.  For all $i,j$, the groups $\Dt_{i}$ and $\Dt_{j}$ commute, and so the group $\Dt=\prod_{i\in I}\Dt_{i}$ is a quotient of the direct product of the $D_{i}$. 
For $a\in \fk^*$ and $i\in I$, let 
\begin{align*}
d_i(a)=\begin{pmatrix}
a & 0 \\0 & a^{-1}
\end{pmatrix}\in D_{i}\sbe \SL_2(\fk)
\end{align*}
and let $\hat{d}_i(a)$ be its image in $\Dt_{i}$.
\newpage
Construct the amalgam of subgroups of $\G^\delta$
\begin{align}\label{eqn:amB}
\amB =\{B_J=\L_J\Dt\mid J\sbne I\}.
\end{align}
Because the group $\Dt_{i}$ either centralizes or normalizes $\L_{j}=\L_{\{j\}}$ for all $j$ we obtain that $\L_J$ is normal in $B_J$. Moreover the action of $\Dt$ on $\L_J$ is induced by the action of $\Dt$ on the $\L_{i}$ so it is determined by the amalgam $\amG$. Since the groups $\L_J$ are perfect, they are contained in $[B_J, B_J]$ and since $B_J/\L_J= \Dt/(\Dt\cap \L_J)$ is abelian, $[B_J, B_J]= \L_J$.

We need to investigate the structure of these groups.
Recall our indexing convention $(i)=I-\{i\}$ for all $i\in I$.  In particular, the maximal groups $B_{(i)}=\langle \L_{(i)},\Dt_i\rangle$ are described by the following lemma.

\ble\label{lem:H}
For any $i$, we have 
$B_{(i)}/H(\G^\delta)\cong (\L_{(i)}\rtimes \Dt_{i})/ H_i$ where  
\begin{align*}
H(\G^\delta)&=\{\hat{d}_1(a)\hat{d}_2(a) \cdots \hat{d}_n(a) \mid a =a^\delta\}\le Z(\G^\delta),\\
H_i(\G^\delta)&=\{(\hat{d}_1(a)\hat{d}_2(a) \cdots \hat{d}_n(a)\hat{d}_i(a)^{-1},\hat{d}_i(a))\mid a=a^\delta\}.
\end{align*}\ele

\bpf
Since the diagram $\Gamma$ is symmetric and $\Dt$ is commutative, we may assume that $i=1$.

Let us consider $d_1(a) \in D_1$ such that $\hat{d}_1(a)$ belongs to $\L_{(1)}\cap \Dt_1$. 
Note that $\hat{d}_1(a)$ commutes with $\L_{j}$ unless $j =1,2$ or $n$. Therefore we need to look at the conjugacy action of $\hat{d}_1(a)$ on $\L_2$ and $\L_n$. Using the definition of $\amG$ we note that $\hat{d}_1(a)$ acts as 
$\begin{pmatrix}
a^{-1} & 0 \\0 & 1
\end{pmatrix}$ on $\L_2$ and as $\left(\begin{array}{cc}1 & 0 \\0 & a^\delta\end{array}\right)$ on $\L_n$ (here we shall write $a^\delta = a^\kaut $ if $\delta =(\kaut,1)$ and $a^\delta =(a^{-1})^\kaut$ if $\delta = (\kaut, \tau)$). In other words, $\hat{d}_1(a)$ acts on $\L_{(1)}$ the same way as the element
\begin{align*} 
d'(a)=
\begin{pmatrix}
a^{-1} &  &\\  & I_{n-2} &  \\ &  & a^\delta
\end{pmatrix}
\end{align*} 
and so, since $d'(a) \in \L_{(1)}$, which is a quotient of $\SL_n(\fk)$, we have $a^\delta =a$ and $d'(a)=(\hat{d}_n(a))^{-1}\cdots (\hat{d}_2(a))^{-1}$. 

More generally, assume $a\in \fk$ is any element satisfying $a^\delta=a$.
This means that the product $\hat{d}(a)=\hat{d}_1(a)({d'}(a))^{-1}= \hat{d}_1(a)\hat{d}_2(a) \cdots \hat{d}_n(a)$ acts trivially on $\L_{(1)}$. Moreover note that the $\hat{d}_i(a)$ commute  and so if $g \in \L_1$, the element $g^{\hat{d}(a)}=(\hat{d}(a))^{-1}g\hat{d}(a)= g^{\hat{d}_1(a)\hat{d}_2(a)\hat{d}_n(a)}$ because the other $\hat{d}(a)$'s commute with $g$. Moreover $\hat{d}_1(a)$, $\hat{d}_2(a)$, $\hat{d}_n(a)$, and $g$ are all in $\L_{\{1,2,n\}}$ and an immediate computation inside this group shows that in fact $g^{\hat{d}(a)}=g$. This shows that $H(\G^\delta)\le Z(\G^\delta)$. 

Now consider the natural homomorphism
$\pi\colon \L_{(1)}\rtimes D_1 \to B_{(1)}/H(\G^\delta)$. 
Clearly $H_1(\G^\delta)\le \ker \pi$.
Now suppose that $(x,y)\in \ker \pi$. Then $y=\hat{d}_1(b)$ for some $b\in \fk$ and $xy=\hat{d}_1(a)\hat{d}_2(a) \cdots \hat{d}_n(a)$ for some $a$ with $a=a^\delta$.
It follows that $x=\hat{d}_1(ab^{-1})\hat{d}_2(a) \cdots \hat{d}_n(a)$ and so $\hat{d}_1(ab^{-1})\in \L_{(1)}$.
From the preceding argument it follows that $(ab^{-1})^\delta=ab^{-1}$ and therefore
$$x\hat{d}_2(b^{-1})\hat{d}_3(b^{-1}) \cdots \hat{d}_n(b^{-1})\in \L_{(1)}\cap H(\G^\delta)
=\{1\}$$
so that $x=\hat{d}_2(b)\hat{d}_3(b)\cdots \hat{d}_n(b)$.
Thus $xy\in H_1(\G^\delta)$.
\epf

From now on, we will let $H(\G^\delta)$ be the group constructed in Lemma~\ref{lem:H}, for any completion $\G^\delta$ of $\amG_{|n-1|}$.

\bpr\label{prop:efficient completion recognition}
Let $\amG^\delta$ be the Curtis-Tits amalgam of type $\tA_{n-1}$ of Definition~\ref{dfn:standard CT structure}. Suppose $\G^\delta$ is a group such that 
\begin{enumerate}
\item $\G^\delta$ contains groups $\L_{i}, \L_{i,j}$ so that the  amalgam $\amL=\{\L_{i}, \L_{i,j}\mid i,j\in I\}$ is isomorphic to $\amG^\delta$
\item $H(\G^\delta)$ is trivial,
\item $\G^\delta$ is the universal completion of the amalgam $\amB$ obtained from $\amL$ as above.		
\end{enumerate}
then  the universal completion $\tG$ of $\amG^\delta$ is an  extension of $\G^\delta$ by $H(\tG)\le Z(\tG)$.
  \epr
\bpf Let  $\tG$ the universal completion of $\amG_{|n-1|}$.
Note that since $\L_i\cong\SL_2(\fk)$, the same is true of the image of $\G_i$ in $\tG$, so that in particular $H_i(\G^\delta)=H_i(\tG)$ for all $i\in I$.

Consider the group $\hG=\tG/H(\tG)$ which is also a completion of $\amG_{|n-1|}$. 
By Lemma~\ref{lem:H} and the observation just made,  $\hG$ is a completion of  $\amB$ and so there is a unique surjective map $\G^\delta\twoheadrightarrow \hG$. Conversely, note that $\L_J$ is isomorphic to the derived subgroup of $B_J$ and so the group $\G^\delta$ contains a copy of the amalgam $\amG_{|n-1|}$. This gives a map $\tG \twoheadrightarrow \G^\delta$. By construction, the map factors through $H(\tG)$. The two maps are inverses to one another since their compositions are the identity on the corresponding amalgams $\amB$ and $\amG_{|n-1|}$.
\epf
In the rest of the paper we will construct a group $\G^\delta$ for any $\delta \in \Aut(\fk)\times \langle \tau \rangle$.
\section{Orientable Curtis-Tits groups}\label{section:orientable}

\subsection{Twisted Laurent polynomial ring \texorpdfstring{$\R_\kaut$}{}, division ring of fractions \texorpdfstring{$\Q_\kaut$}{}, and linear groups}\label{subsec:twisted Laurent polynomial ring}

Recall that $\fk$ is a commutative field of order at least $4$ and that $\kaut\in \Aut(\fk)$.
If $\kaut$ has finite order $s$, let $T=t^s$ and let $\A=\fk[T,T^{-1}]\le \R_\kaut$ be the ring of Laurent polynomials in the commuting variable $T$ with coefficients in the commutative field $\fk$. Moreover, let $\F=\fk(T)$.

As $\fk\{t\}=\fk[t,\kaut^{-1}]$, in the notation of~\cite{Jac1943a}, is a (non-commutative) principal ideal domain, it is in particular a left and right Ore ring, and so possesses a division ring of fractions, which we shall denote $\Q_\kaut$ (see also~\cite{Coh1971}). 
Naturally, $\R_\kaut\le \Q_\kaut$. Also, for finite $s$, identify $\F$ with the subfield of $\Q_\kaut$ generated by $\fk$ and $T$.
Let $\V$ be a left $\Q_\kaut$-vector space of dimension $n$ and  $\M\le \V$ a free $\R_\kaut$-submodule of rank $n$, so that $\Q_\kaut \M=\V$.
The group of all $\Q_\kaut$- (resp. $\R_\kaut$-) linear invertible transformations of $\V$ (resp. $\M$) is denoted $\GL_{\Q_\kaut}(\V)$ (resp. $\GL_{\R_\kaut}(\M)$).

We fix an ordered reference $\Q_\kaut$-basis $\cE=\{e_1,e_2,\ldots,e_n\}$ of $\V$ that is also an $\R_\kaut$-basis for $\M$. 
We will represent an element $x=\sum_{i=1}^n x_ie_i\in \V$ as a row vector $(x_1,\ldots,x_n)$.
Representation of $\Q_\kaut$-linear endomorphisms of $\V$  as matrices w.r.t.~the basis $\cE$ by matrix multiplication on the right yields the usual
identification: $\End_{\Q_\kaut}(\V) \to M_n(\Q_\kaut)$.
The images  of $\GL_{\Q_\kaut}(\V)$ and $\GL_{\R_\kaut}(\M)$ under this identification will be denoted
  $\GL_n({\Q_\kaut})$ and $\GL_n({\R_\kaut})$ respectively. 
The inclusion $\cE\sbe \M\sbe \V$ induces the  inclusions
 $\GL_{\R_\kaut}(\M)\le \GL_{\Q_\kaut}(\V)$ and $\GL_n({\R_\kaut})\le \GL_n({\Q_\kaut})$.

The Dieudonn\'e determinant (see~\cite{Die1943a}) is the unique non-trivial group homomorphism 
\begin{align}
\Det\colon \GL_n({\Q_\kaut})\to {\Q_\kaut^*}/[{\Q_\kaut^*},{\Q_\kaut^*}]
\end{align}
which is trivial on transvections, and induces the canonical homomorphism ${\Q_\kaut^*}\to~{\Q_\kaut^*}/[{\Q_\kaut^*},{\Q_\kaut^*}]$ on diagonal matrices having exactly one non-identity entry.
Here $[{\Q_\kaut^*},{\Q_\kaut^*}]$ denotes the commutator subgroup of the multiplicative group ${\Q_\kaut^*}$.
If ${\Q_\kaut}$ is commutative $\Det$ is just the ordinary determinant.

We let $\SL_n({\Q_\kaut})$ (resp. $\SL_n({\R_\kaut})$, $\SL_{\R_\kaut}(\M)$, $\SL_{\Q_\kaut}(\V)$) be the kernel of 
 $\Det$ restricted to $\GL_n({\Q_\kaut})$ (resp. $\GL_n({\R_\kaut})$, $\GL_{\R_\kaut}(\M)$, $\GL_{\Q_\kaut}(\V)$).

\bde\label{dfn:norms}
Recall that $\fk_\kaut$ is the fixed field of $\kaut$ in $\fk$.
Assume that $\kaut$ has finite order $s$.
We denote the image of the norm map $N_{\fk/\fk_\kaut}\colon b\mapsto \prod_{i=0}^{s-1}b^{\kaut^i}$ by $\n_\kaut\le \fk_\kaut^*$.
This extends to a norm map $N_{{\R_\kaut^*}/\A^*}\colon bt^k\mapsto N_{\fk/\fk_\kaut}(b)((-1)^{s-1}T)^k$, where $T=t^s$. 
Note that this is the restriction of the standard reduced norm for the cyclic algebra $\R_{\kaut}$ over $\fk(T)$. More precisely, $(-1)^{s-1}T$ is the determinant of the image of $t$ under the splitting morphism from $\R_\kaut$ to $M_s(\overline{\fk(T)})$.
\ede
\ble\label{lem:norms}
We have
\begin{enumerate}
\item ${\R_\kaut^*}=\{bt^l\mid b\in \fk, l\in \ZZ\}$,
\item $[{\R_\kaut^*},{\R_\kaut^*}]=\langle b^{\kaut^l}b^{-1}\mid b\in \fk, l\in \ZZ\rangle=\{b^\kaut b^{-1}\mid b\in \fk\}$
\item $N_{{\R_\kaut^*}/\A^*}$ induces a surjective homomorphism 
\begin{align*}
{\R_\kaut^*}/[{\R_\kaut^*},{\R_\kaut^*}]\to \{n((-1)^{s-1}T)^l\mid n\in \n_\kaut,l\in \ZZ\},
\end{align*} which is an isomorphism provided $\fk/\fk_\kaut$ is a separable (hence cyclic Galois) extension.
\end{enumerate}
\ele 
\bpf
(a) ``$\spe$'' is clear. For the converse note that if $f\in {\R_\kaut}$ has at least two terms, then so does any multiple of $f$ and so $f$ cannot be a unit.
(b) The first equality follows from (a) by direct computation.
For the second equality, note that since $\fk$ is commutative, for $l\ge 1$, 
\begin{align*}
b^{\kaut^{l}}b^{-1} & = \prod_{i=0}^{l-1}  (b^{\kaut^{i}})^{\kaut} (b^{\kaut^{i}})^{-1}
= (\prod_{i=0}^{l-1}  b^{\kaut^{i}})^{\kaut} (\prod_{i=0}^{l-1}  b^{\kaut^{i}})^{-1}.
\end{align*}
(c) Since conjugate elements have the same norm, this map is a homomorphism. Surjectivity is obvious. Injectivity follows from Hilbert's 90th theorem. 
\epf

Let $\Z_n({\R_\kaut})=\Z(\GL_n({\R_\kaut}))$.
Define $\PGL_n({\R_\kaut})=\GL_n({\R_\kaut})/\Z_n({\R_\kaut})$ and $\PSL_n({\R_\kaut})=\SL_n({\R_\kaut})/(\Z_n({\R_\kaut})\cap\SL_n({\R_\kaut}))$.
We shall interpret $\PSL_n({\R_\kaut})$ as a subgroup of $\PGL_n({\R_\kaut})$ via $\PSL_n({\R_\kaut})\cong \SL_n({\R_\kaut})\cdot \Z_n({\R_\kaut})/\Z_n({\R_\kaut})$.

\bpr\label{prop:PGL:PSL}
Let $\fk/\fk_\kaut$ be a cyclic Galois extension. Then, we have 
\begin{align*}
|\PGL_n({\R_\kaut})\colon \PSL_n({\R_\kaut})|=sn|\n_\kaut\colon (\fk_\kaut^*)^{sn}|.
\end{align*}
\epr
\bpf
We shall make use of the fact that 
\begin{align*}
|\PGL_n({\R_\kaut})\colon \PSL_n({\R_\kaut})|
=
|\GL_n({\R_\kaut})\colon \SL_n({\R_\kaut})\Z_n(\R_\kaut)|.
\end{align*}
Consider the composition $\chi$ of surjective homomorphisms (compare Lemma~\ref{lem:norms}):
\begin{align*}
\GL_n({\R_\kaut})\stackrel{\Det}{\to}{\R_\kaut^*}/[{\R_\kaut^*},{\R_\kaut^*}]\stackrel{N_{{\R_\kaut^*}\! /\A^*}}{\to}\!\{n((-1)^{s-1}T)^k\mid n\in \n_\kaut, k\in \ZZ\}\cong \n_\kaut\times\ZZ.
\end{align*}
We claim that 
\begin{align*}
\Z_n({\R_\kaut})=\{b t^{sl}I_n\mid b\in \fk_\kaut, l\in \ZZ\},
\end{align*}
where $I_n$ denote the $n\times n$ identity matrix. The inclusion $\spe$ is clear since $bt^{sl}\in \Z({\R_\kaut^*})$.
Conversely, by considering commutators with permutation matrices, it follows that a central element in $\GL_n({\R_\kaut})$ must be scalar. It then follows that the scalar must belong to the center $\Z({\R_\kaut^*})$.
Now $\chi(bt^{sl}I_n)=N_{{\R_\kaut^*}/\A^*}(b^nt^{snl}[{\R_\kaut^*},{\R_\kaut^*}])=b^{sn}((-1)^{s-1}T)^{snl}$ and since $\{n((-1)^{s-1}T)^k\mid n\in \n_\kaut, k\in \ZZ\}\cong \n_\kaut\times\ZZ$ we see that 
\begin{align*}
\GL_n({\R_\kaut})/\SL_n({\R_\kaut})\cdot \Z_n({\R_\kaut})\cong \n_\kaut/(\fk_\kaut^*)^{sn} \times\ZZ/{sn}\ZZ.
\end{align*}
\epf
\subsection{A realization of \texorpdfstring{$\amG^\kaut$}{} inside \texorpdfstring{$\SL_n({\R_\kaut})$}{}.}\label{subsec:realization of G alpha}
At the very end of~\cite{Ti1992} it is claimed that a  Kac-Moody group $\G^\kaut$ that is a completion of $\amG^\kaut$ can be obtained as a subgroup inside $\PGL_n({\R_\kaut})$. We shall now proceed to give an explicit description of the amalgam inside $\SL_n(\R_\kaut)$. Since the amalgam does not intersect the center, this gives rise to a realization inside $\PSL_n(\R_\kaut)$, which, in turn, via Proposition~\ref{prop:PGL:PSL} can be viewed as a subgroup of index $sn|\n_\kaut\colon (\fk_\kaut^*)^{sn}|$ inside $\PGL_n(\R_\kaut)$.

In order exhibit this amalgam, we first define the following injective homomorphisms $\phi_i\colon \SL_2(\fk)\into \GL_n(\R_\kaut)$.
For $i=1,\ldots,n-1$ we take 
\begin{align*}
\phi_i\colon A\mapsto \left(\begin{array}{@{}ccc@{}}
 I_{i-1} &   &                    \\
      &A &                    \\
      &    & I_{n-i-1}    \\
\end{array}\right).
\end{align*}
Moreover, we define
\begin{align*}
\phi_{n}\colon \left(\begin{array}{@{}cc@{}} a & b \\ c & d \end{array}\right) & \mapsto \left(\begin{array}{@{}ccc@{}}
  d^{\kaut^{-1}}         &                  &  t^{-1}c\\
                          & I_{n-2}     &  \\
 bt                  &                 & a \\
 \end{array}\right).
 \end{align*}
Now, for every $i\in I$, let $\L_{i}=\im \phi_i$ and 
 $\L_{i,j}=\langle \L_i,\L_j\rangle\le {\GL_n(\R_\kaut)}$.
Consider the amalgam $\amL^\kaut({\R_\kaut})=\amL^\kaut=\{\L_{i},\L_{i,j}\mid i,j\in I\}$ of subgroups of $\GL_{n}({\R_\kaut})$.
Here the connecting maps $\varphi_{i,j}$ of $\amL^\kaut$ are the natural inclusion maps of subgroups of $\GL_n({\R_\kaut})$.

\bpr\label{prop:oriented CT amalgam}
We have an isomorphism of amalgams
 $\amL^\kaut\cong \amG^\kaut$.
Hence, $\G^\kaut=\langle \amL^\kaut\rangle$ is a non-trivial completion of $\amG^\kaut$ inside $\SL_n({\R_\kaut})$.
 \epr
\bpf
Consider the following matrix:
\begin{align}\label{eqn:\Cyc_R}
\Cyc=\Cyc_{{\R_\kaut},n}=\left(\begin{array}{@{}c|c@{}}
0 & I_{n-1} \\
\hline
t & 0 \\
\end{array}\right).
\end{align}
We now define the automorphism $\Phi=\Phi_{\R_\kaut}$ of $\GL_n({\R_\kaut})$ given by $X\mapsto \Cyc^{-1}X\Cyc$.
One verifies that, for $i=1,\ldots,n$ we have
 $\phi_i=\Phi^{i-1}\after\phi_1$.
In particular $\phi_{n}$ is an isomorphism.
We now turn to the rank $2$ groups.
For distinct $i,j\in \{1,2,\ldots,n\}$, let $\phi_{i,j}$ be the canonical isomorphism between 
 $\G_{i,j}=\langle \G_{i},\G_{j}\rangle$ and $\L_{i,j}=\langle \L_{i},\L_{j}\rangle$ induced by $\phi_i$ and $\phi_{j}$. Note that this implies that $\phi_{i,i+1}=\Phi^{i-1}\after\phi_{1,2}$.

We claim that the collection $\phi=\{\phi_i,\phi_{i,j}\mid i,j\in I\}$ is the required isomorphism between $\amG^\kaut$ and $\amL^\kaut$.
This is completely straightforward except for the maps
 $\phi_1$, $\phi_{n,1}$.
Note that 
\begin{align*}
\phi_{n,1}\colon
\left(\begin{array}{@{}ccc@{}}
a & b & c \\
d & e & f\\
g & h & i\\
\end{array}\right)\mapsto&
\left(\begin{array}{@{}cc|c|c@{}}
t^{-1}e t & t^{-1}f t &  & t^{-1}d \\
t^{-1}h t & t^{-1}i t &  & t^{-1}g \\
\hline
             &              & I_{n-3} &  \\
\hline
bt & ct &     & a\\
\end{array}\right)
\\
=&\left(\begin{array}{@{}cc|c|c@{}}
e^{\kaut^{-1}}  & f^{\kaut^{-1}}  &  & t^{-1}d \\
h^{\kaut^{-1}} & i^{\kaut^{-1}}  &  & t^{-1}g \\
\hline
             &              & I_{n-3} &  \\
\hline
bt & ct &     & a\\
\end{array}\right).
\end{align*}
Thus we have
$$\phi_{i,j}\after \psi_{i,j} =\varphi_{i,j}\after  \phi_i,$$
for all $i,j\in I$.

Since all $\L_{i}$ are conjugates of $\L_1$, which clearly lies in $\SL_n({\R_\kaut})$ and the Dieudonn\'e determinant is a homomorphism to the abelian group ${\R_\kaut^*}/[{\R_\kaut^*},{\R_\kaut^*}]$, the second claim follows.
\epf
\subsection{The twin-building of type \texorpdfstring{$\tA_{n-1}$}{} over \texorpdfstring{${\R_\kaut}$}{}}\label{subsec:building over R}
We take the excellent and succinct description from~\cite{AbVa2001} and adapt it to the non-commutative setting we need.
Let $\valuation_+,\valuation_-\colon {\Q_\kaut}\to \ZZ$ be the non-commutative discrete valuations determined by $\valuation_+(\fk^*)=\valuation_-(\fk^*)=0$ and  $\valuation_+(t)=\valuation_-(t^{-1})=1$, and let $\cO_\vep=\{\lambda\in {\Q_\kaut}\mid \valuation_\vep(\lambda)\ge 0\}$ ($\vep=+,-$) be the corresponding valuation ring.

An $\cO_\vep$-lattice is a free left $\cO_\vep$ module $\Lat\le \V$ with ${\Q_\kaut}\Lat=\V$.
Such lattices are of the form
\begin{align*}
\Lat=\bigoplus_{i=1}^n \cO_{\vep} a_i,
\end{align*}
where $\{a_1,a_2,\ldots,a_n\}$ is a ${\Q_\kaut}$-basis for $\V$.
We call $\{a_1,a_2,\ldots,a_n\}$ a {\dfn lattice basis} for $\Lat$.

A chain  $\cdots\sbne \Lat_i\sbne \Lat_{i+1}\sbne\cdots$ of $\cO_\vep$-lattices is called {\dfn admissible} if it is invariant under multiplication by integral powers of $t$.
The admissible chain generated by the lattice $\Lat$ will be denoted $[\Lat]$.

For $\vep=+,-$, we now describe an incidence geometry $\geomI_\vep$.
The {\dfn objects} of $\geomI_\vep$ are the minimal admissible chains of $\cO_\vep$-lattices; these are of the form $\LatU=[\Lat]$ for some lattice $\Lat$.
Call two objects $\LatU$ and $\LatU'$ {\em incident} if $\LatU\cup \LatU'$ is admissible.
Naturally, a flag is given by a set $\{\LatU_1,\ldots,\LatU_r\}$ of objects such that $\LatU_1\cup\cdots \cup \LatU_r$ is admissible.
The {\em chambers} of $\geomI_\vep$ are maximal flags.
Following loc.~cit.~we associate the following to any ordered ${\Q_\kaut}$-basis $(a_1,\ldots,a_n)$ of $\V$ and $j\in \{0,1,\ldots,n-1\}$:
\begin{align*}
\Lat_\vep^j(a_1,\ldots,a_n)&:=\langle ta_1,\ldots,ta_j,a_{j+1},\ldots,a_n\rangle_{\cO_\vep},\\
\LatU_\vep^j(a_1,\ldots,a_n)&:=[\Lat_\vep^j],\\
c_\vep(a_1,\ldots,a_n) & := \{\LatU_\vep^0,\ldots,\LatU_\vep^{n-1}\}.
\end{align*}
The latter is called the chamber with {\em ordered chain basis} $(a_1,\ldots,a_n)$.

The geometry $\geomI_\vep$ has type set $\{0,1,\ldots,n-1\}$.
The {\em type} function is given by 
$\typ_{\vep}([\Lat_\vep^0(g(e_1),\ldots,g(e_n))])={\vep \nu_\vep(\Det(g))\mod n }$ for all $g\in \GL_{\Q_\kaut}(\V)$, where $\Det$ denotes the Dieudonn\'e determinant.
In particular,  $\typ_\vep(\LatU_\vep^j(e_1,\ldots,e_n))=j$, for $j=0,1,\ldots,n-1$.

Let $\De_\vep$ be the chamber system of $\geomI_\vep$ in which two chambers $c_\vep$ and $d_\vep$ are {\dfn $i$-adjacent}, written $c_\vep\sim_i d_\vep$, if their objects of type $j\ne i$ are equal.

Given a ${\Q_\kaut}$-basis $\{a_1,\ldots,a_n\}$ for $\V$, we define the subsystem 
\begin{align*}
\Sigma_\vep(a_1,\ldots,a_n):=\{c_\vep(t^{m_1}a_1,\ldots,t^{m_n} a_n)\mid m_1,\ldots,m_n\in \ZZ\}.
\end{align*}

It can be proved (see e.g.~\cite[\S 9.2]{Ro1989a}) that $\De$ with given adjacency relations forms a building of affine type $\tA_{n-1}(\fk)$ and that the collection $$\cA_\vep=\{\Sigma_\vep(a_1,\ldots,a_n)\mid (a_1,\ldots, a_n)  \mbox{ is a ${\Q_\kaut}$-basis for $\V$}\}$$ is a system of apartments for $\De_\vep$.

We now define a symmetric {\dfn opposition relation}  $\opp\sbe \De_+\times \De_-\cup \De_-\times\De_+$ by declaring 
 $c_+\opp c_-$ if and only if 
  $c_\vep=c_\vep(a_1,\ldots,a_n)$ ($\vep=+,-$) for some ${\R_\kaut}$-basis $\{a_1,\ldots,a_n\}$ for $\M$.
Moreover, two objects are declared {\em opposite} if they belong to opposite chambers and have the same type.

The proof given in~\cite[\S 4]{AbVa2001}, which is given in the context where ${\Q_\kaut}$ is commutative, can be applied almost verbatim to prove the following.
\bpr\label{prop:De twin building of type tA_{n-1}}
$(\De_+,\De_-,\opp)$ is a twin-building of type $\tA_{n-1}(\fk)$ with system of twin-apartments 
\begin{align*}
\cA_{\opp}&=\{(\Sigma_\vep(a_1,\ldots,a_n)\colon \vep=\pm)\mid   (a_1,\ldots,a_n)\mbox{ is an ${\R_\kaut}$-basis for $\M$}\}.
\end{align*}
\epr
\bre
The group $\GL_{\R_\kaut}(\M)$ is a group of sign-preserving automorphisms of $(\De_+,\De_-,\opp)$, which does not preserve types.
\ere

\ble\label{lem:SL flag transitive}
The group $\SL_{\R_\kaut}(\M)$ of type preserving automorphisms of the twin-building $(\De_+,\De_-,\opp)$ acts transitively on pairs of opposite chambers.
\ele
\bpf
For $\vep=\pm$, $\SL_{\R_\kaut}(\M)$ is a group of permutations of the collection of $\cO_\vep$-lattices that preserve containment and types.
Suppose $(c_+,c_-)$ and $(d_+,d_-)$ are pairs of opposite chambers. 
Without loss of generality assume that $c_\vep=c_\vep(e_1,\ldots,e_n)$ and $d_\vep=c_\vep(b_1,\ldots,b_n)$ for a suitable ordered ${\R_\kaut}$-basis 
$(b_1,\ldots,b_n)$ for $\M$ and $\vep=+,-$.
Then there is $g\in \GL_{\R_\kaut}(\M)$ with $g(e_i)=b_i$ for $i=1,2,\ldots,n$.
Let $\Det(g)$ be represented by $a t^m$ in ${\R_\kaut^*}/[{\R_\kaut^*},{\R_\kaut^*}]$ for some $a\in \fk$ and $m\in \ZZ$.
Since $\LatU_\vep^0(e_1,\ldots,e_n)$ and $\LatU_\vep^0(b_1,\ldots,b_n)$ have type $0$ apparently $\vep\valuation_\vep(\Det(g))=0\mod n$ so that $m=nl$ for some $l\in \ZZ$.
This means that $g'\in \GL_{\R_\kaut}(\M)$ given by 
 $g'(e_1)=a^{-1}t^{-l}b_1$, $g'(e_i)=t^{-l}b_i$ ($i=2,3, \ldots, n$) also satisfies $g'(c_+,c_-)=(d_+,d_-)$.
Also, $\Det(g')=\Det(g)\cdot a^{-1}t^{-m}\in [{\R_\kaut^*},{\R_\kaut^*}]$, so that $g'\in \SL_{\R_\kaut}(\M)$.
\epf

Let $\GD$ (resp.\ $\Dt$) be the maximal split torus in $\GL_n({\R_\kaut})$ (resp.\ $\SL_n({\R_\kaut})$) stabilizing the pair of opposite chambers $(c_+,c_-)$, where $c_\vep=c_\vep(e_1,\ldots,e_n)$.
The group $\Dt$ is generated by the images $\Dt_i$ ($i\in I$) of $D_i$ and so it appear in the definition of $\amB$ as in~\eqref{eqn:amB} when we apply Proposition~\ref{prop:efficient completion recognition}.

\ble\label{lem:stabilizer of c+ and c-}
Let $c_\vep=c_\vep(e_1,\ldots,e_n)$ for $\vep=\pm$.
\begin{enumerate}
 \item
 The stabilizer $\Dt$ of $(c_+,c_-)$ in $\SL_n({\R_\kaut})$ is the subgroup of diagonal matrices of Dieudonn\'e determinant $1$ and coefficients in $\fk$.
\item The stabilizer $\GD$ of $(c_+,c_-)$ in $\GL_n({\R_\kaut})$ is the subgroup generated by diagonal matrices in $\fk^*$ and 
 scalar matrices with coefficients in ${\R_\kaut^*}$.

\end{enumerate}
\ele
\bpf

(a) Let $g\in\SL_{\R_\kaut}(\M)$ preserve $c_+$ and $c_-$. Then, $g$ stabilizes the objects  
$\LatU^0_\vep(e_1,\ldots,e_n)$, for $\vep=\pm1$.
Since $\Det(g)=1$, $g$ preserves the intersection 
$\Lat^0_+(e_1,\ldots,e_n)\cap \Lat^0_-(e_1,\ldots,e_n)$ and so $g\in \GL_n(\fk)$.
Now, $g$ preserves two opposite chambers in the $0$-residue on $c_+$, which is the spherical building $\Lat^0_+/t\Lat^0_+$ of type $A_{n-1}(\fk)$. This shows that $\Dt$ is contained in the group of diagonal matrices in $\GL_n(\fk)$ with Dieudonn\'e determinant $1$.
Conversely, note that the images $\Dt_i$ of the $D_{i}$ ($i=1,2,\ldots,n$) generate $\Dt$.
Now the description of $D_n$ together with Lemma~\ref{lem:norms} shows that $\Det(\Dt)= [{\R_\kaut^*},{\R_\kaut^*}]$.

(b) Let $g'\in \GL_{\R_\kaut}(\M)$ preserve $c_+$ and $c_-$
Then, $\Det(g)=a t^{ln}/[{\R_\kaut^*},{\R_\kaut^*}]$ for some $a\in \fk^*$ and $l\in \ZZ$ since $g'$ preserves the type of the $0$-object on $c_+$. Define $d\in \GL_{\R_\kaut}(\M)$ by $d(e_1)=a^{-1}t^{-l}e_1$, 
 and $d(e_i)=t^{-l}e_i$.
Then, $\Det(g)=\Det(d)\Det(g')=1/[{\R_\kaut^*},{\R_\kaut^*}]$ so $g\in \SL_{\R_\kaut}(\M)$ and the result follows from (a).
\epf

\bpf (of Theorem~\ref{mainthm:orientable CT group})
By Proposition~\ref{prop:De twin building of type tA_{n-1}} $\De$ is a twin-building with diagram $\tA_{n-1}$, where $n\ge 4$.
In particular, $\De$ satisfies condition (co) of~\cite{MuRo1995b}. 
By Lemma~\ref{lem:SL flag transitive}, $\SL_{\R_\kaut}(\M)$ is an automorphism group of ${\De}$ that is transitive on pairs of opposite chambers.
Define the  amalgam 
$\amB_2=\{B_{i},B_{ij}\mid i,j\in I\}$ of  Levi-components of rank $1$ and $2$ and the amalgam
$\amB=\{B_J=\langle B_{i}\mid i\in J\rangle\mid J\sbne I\}$.
Then, by the twin-building version of the Curtis-Tits theorem~\cite{AbMu1997} the automorphism group $\SL_{\R_\kaut}(\M)$ of ${\De}$ is the universal completion of 
$\amB_2$ and, a forteriori $\SL_{\R_\kaut}(\M)$ is the universal completion of the amalgam $\amB$.
Now consider the amalgam $\amL^\kaut$. One verifies easily that, for each $i,j\in I$, $\SL_2(\fk)\cong \L_{i}\le B_{i}$
 and $\SL_3(\fk)\cong \L_{ij}\le B_{ij}$, when $\{i,j\}$ is an edge of the diagram.
 In fact for any $J\sbne I$, we have $B_J=\L_J\Dt$; this follows for instance by considering the transitive action of both groups on the pair of opposite residues of type $J$ on $(c_+,c_-)$.
This means that $\amB$ is defined as in~\eqref{eqn:amB} and so, in view of Proposition~\ref{prop:efficient completion recognition}, it suffices to show that $H(\L)=1$.
This follows by noting that if $a=a^\kaut$, then taking the product over all $\phi_i$ images of the matrix
 \begin{align*}
 \begin{pmatrix}
 a & 0\\
 0 & a^{-1}
 \end{pmatrix}
 \end{align*}
we obtain the identity of $\SL_n({\R_\kaut})$.
\epf

\bre
Note that this construction is in particular valid if $\kaut=\id$ and the classical definition of the building over commuting Laurent polynomials follows.
Thus, in the above, we can replace the skew Laurent polynomial ring ${\R_\kaut}$ and its division ring of fractions $\Q_{\kaut}$ by the Laurent polynomial ring $\A$ and its field of fractions $\F$ (see the definitions at the beginning of Subsection~\ref{subsec:twisted Laurent polynomial ring}).
Note that in that case, where $\kaut=\id$, a slightly weaker statement in the vein of Theorem~\ref{mainthm:orientable CT group} can be deduced from~\cite{Cap2007}.
\ere

\section{The non-orientable Curtis-Tits group \texorpdfstring{$\G^\delta$}{}}\label{section:non-orientable}
We adopt the notation of Section~\ref{section:ct-groups}~and~\ref{section:orientable}. We assume that $\delta=\kaut\tau$ has finite order $s$.
As in Section~\ref{section:orientable}, ${\R_{\kaut^2}}=\fk\{t,t^{-1}\}$ denotes the ring of not necessarily commuting Laurent polynomials with coefficients in the field $\fk$.
Here, for $b\in \fk$, we have $tbt^{-1}=b^{\kaut^2}$.

Let $I=\{1,2,\ldots,n\}$ and let $\tI=\{1,2,\ldots,2n\}$.
As before let $\V$ be a left ${\Q_{\kaut^2}}$-vector space of dimension $2n$, where $n\ge 4$, with (ordered) basis $\cE=\{e_{1}, \dots, e_{n}, f_{1}=e_{n+1},\ldots, f_{n}=e_{2n}\}$.
The vector $x=\sum_{i=1}^{2n} x_i e_i$ will be represented as the row vector $(x_1,\ldots,x_{2n})$.
Let $\M$ be the free ${\R_{\kaut^2}}$-module spanned by this basis. 
As in Section~\ref{section:orientable} we identify 
 $\End_{\R_{\kaut^2}}(\M)$ with $M_n({\R_{\kaut^2}})$ via the right action on $\V$.
Furthermore we let $\G=\SL_{\R_{\kaut^2}}(\M)$.

In this subsection we introduce a sesquilinear form $\beta$ on $\V$ and an involution $\theta$ of $\G$ such that the fixed group $\G^\theta$ is precisely the group of symmetries of $\beta$ in $\G$.
In Subsection~\ref{subsection:flag transitive} we will prove that $\G^\theta$ is flag-transitive on a geometry $\De^\theta$. In Subsection~\ref{sec:simple-connectedness} we prove that the
geometry $\De^\theta$ is connected and simply connected which 
 by Tits' Lemma implies that the group $\G^\theta$ is the universal completion of the amalgam of maximal parabolics. 
We then apply Proposition~\ref{prop:efficient completion recognition}

\subsection{\texorpdfstring{$\sigma$}{}-sesquilinear forms on \texorpdfstring{$\V$}{}}\label{subsec:beta}

Let $\sigma$ be an anti-automorphism of ${\Q_{\kaut^2}}$ that interchanges $t$ and $t^{-1}$. Thus $\sigma^2$ fixes $t$, but may act as a non-trivial automorphism of $\fk$.

We wish to define a $\sigma$-sesquilinear form $\beta$ on $\V$. This is a function $\beta\colon \V\times \V\to {\Q_{\kaut^2}}$ satisfying 

\begin{align*}
&\beta(\lambda u,\mu v) =\lambda \beta(u,v)\mu^\sigma,\\
&\beta(u_1+u_2, v)=\beta(u_1,v)+\beta(u_2,v), &\\
&\beta(u, v_1+v_2)=\beta(u,v_1)+\beta(u,v_2), &
\end{align*}
for all $u,,u_1,v,v_1,v_2,v\in \V$ and $\lambda,\mu\in {\Q_{\kaut^2}}$

Note that $\beta$ is uniquely determined by the Gram matrix $B=(\beta(e_i,e_j))_{i,j=1}^{2n}$ of $\cE$ with respect to $\beta$. We shall assume that $\beta$ is non-degenerate, that is, $B$ is invertible.

More concretely, 
\begin{align}
\beta(x,y)
=(x_1,\ldots,x_{2n})B\,{}^t\!(y_1,\ldots,y_{2n})^\sigma=\sum_{i,j=1}^{2n} x_i b_{i,j}y_j^\sigma.
\end{align}
\bde\label{dfn:adjoint map}
The {\dfn right adjoint} of a transformation $g\in \GL(\V)$, is the transformation $g^\adj\in\GL(\V)$ such that
\begin{align}\label{eqn:adjoint}
\beta(g(u),v)&=\beta(u,g^\adj(v)) &\mbox{ for all }u,v\in \V.
\end{align} 
The {\dfn inverse adjoint} of a transformation $g\in \GL(\V)$, is the transformation $g^*\in\GL(\V)$ such that
 $\beta(g(u),g^* (v))=\beta(u,v)$ for all $u,v\in \V$.
Clearly, $g^*=(g^{-1})^\adj$.
\ede

\ble 
\begin{enumerate}
\item For any two matrices of compatible dimension $X$ and $Y$, we have 
\begin{align*}
{}^t(XY)^\sigma={}^tY^\sigma \cdot {}^tX^\sigma & \mbox{ and }\\
{}^t({}^tX^\sigma)^\sigma = X^{\sigma^2}.
\end{align*}

\item The map $\GL(\V)\to \GL(\V)$, $x\mapsto x^\adj$ is an anti-isomorphism, which via the right action on $\V$ corresponds to the anti-isomorphism $M_{2n}({\Q_{\kaut^2}})\to M_{2n}({\Q_{\kaut^2}})$  given by 
\begin{align*}
X\mapsto X^\adj= {}^tB^{\sigma^{-1}}\ {}^tX^{\sigma^{-1}}\ {}^t(B^{-1})^{\sigma^{-1}}.
\end{align*}
\item The map $\GL(\V)\to \GL(\V)$ given by $x\mapsto x^*$ is an automorphism, corresponding via the  right action on $\V$ to the automorphism of $M_{2n}({\Q_{\kaut^2}})$ given by 
\begin{align}\label{eq:adj mat c}
X\mapsto X^*= {}^tB^{\sigma^{-1}}\ {}^t(X^{-1})^{\sigma^{-1}}\ {}^t(B^{-1})^{\sigma^{-1}}.
\end{align}
\end{enumerate}
\ele
\bpf
(a) Suppose $X=(x_{i,j})$ and $Y=(y_{j,k})$.
Then the $ki$-entry on both sides is 
$\sum_{j} y_{jk}^\sigma x_{ij}^\sigma$. The second equality is clear.

(b) Since $\beta$ is non-degenerate, $x$ uniquely determines $x^\adj$ via the equality~\eqref{eqn:adjoint} and the property $(xy)^\adj=y^\adj\cdot x^\adj$ follows easily.
As for the matrix identity,
let $u=(u_1,\ldots,u_{2n}), v=(v_1,\ldots,v_{2n})\in \V$.
Suppose $x^\adj$ is represented by the matrix $Y$.
Then, apparently
\begin{align*}
u\, X B\, {}^tv^\sigma
=\beta(x(u),v)
=\beta(u,x^\adj(v))
=u\,B\, {}^t(vY)^\sigma.
\end{align*}
Since $u$ and $v$ are arbitrary, using (a) we find that 
\begin{align*}
XB=B\ {}^tY^\sigma
\end{align*}
and so we find that 
\begin{align*}
Y={}^t(B^{-1} X B)^{\sigma^{-1}} = {}^tB^{\sigma^{-1}}\ {}^tX^{\sigma^{-1}}\ {}^t(B^{-1})^{\sigma^{-1}}.
\end{align*}
Claim (c) follows from (b) noting that $x^*=(x^{-1})^\adj$.
\epf

\bde\label{dfn:theta}
For $B\in \GL_{2n}(\sA)$, 
we define an automorphism $\theta\colon \G \mapsto \G$ by $x\mapsto x^*$. If $x$ corresponds to $X$ under the identification $\G=\SL_{2n}({\R_{\kaut^2}})\le \GL_{2n}({\Q_{\kaut^2}})$, then, $\theta$ is given by 
\beq \label{def:theta}
X\mapsto {}^tB^{\sigma^{-1}}\ {}^tX^{-\sigma^{-1}}\ {}^tB^{-\sigma^{-1}}.
\eeq
Note that with this choice of $B$, $X^\theta$ does belong to $\SL_{2n}({\R_{\kaut^2}})$.
Occasionally we shall write $\theta=\theta_{\R_{\kaut^2}}\in \Aut(\SL_{\R_{\kaut^2}}(\M))$ to distinguish it from $\theta_\delta\in \Aut(\SL_\A(\M))$.
\ede
\bde\label{dfn:SU}  
\begin{align}\label{eqn:SU}
\GU_{\R_{\kaut^2}}(\M,\beta)&:=\{ g \in \GL_{\R_{\kaut^2}}(\M)| \forall x,y \in \M,  \beta(gx,{g}y)=\beta(x,y)\}.\\
\SU_{\R_{\kaut^2}}(\M,\beta)&:=\GU_{\R_{\kaut^2}}(\M,\beta)\cap \SL_{\R_{\kaut^2}}(\M).\nonumber
\end{align}
We let $\GU_n({\R_{\kaut^2}})$ and $\SU_n({\R_{\kaut^2}})$ denote the subgroups of $\GL_n({\R_{\kaut^2}})$ corresponding to $\GU_{\R_{\kaut^2}}(\M,\beta)$ and $\SU_{\R_{\kaut^2}}(\M,\beta)$ respectively via its right action on $\V$.
\ede
\bco\label{cor:SU}
The unitary group $\SU_{\R_{\kaut^2}}(\M,\beta)$ is the fixed group $\G^\theta=\{x\in \G\mid x^\theta=x\}$.
\eco

\subsection{The amalgam \texorpdfstring{$\amL^\delta$}{}}\label{subsec:amL^delta}
We shall continue the terminology from Subsection~\ref{subsec:beta} with the following choices for $\sigma$ and $B$.
As in the Introduction, let $\sigma$ be the anti-automorphism of ${\Q_{\kaut^2}}$ that interchanges $t$ and $t^{-1}$ and acts as $\kaut^{-1}$ on $\fk$, and let  
\begin{align}
B=(\beta(e_i,e_j))=\left(  \begin{array}{c|c} 0_n &  I_n \\ \hline  tI_n & 0_n \end{array}\right)\in \GL_{2n}({\R_{\kaut^2}}).
\end{align}
We first note that 
\begin{align}
{}^tB^{-\sigma^{-1}}&=B,\\
X^\theta&=B^{-1}\, {}^tX^{-\sigma^{-1}}\, B &&\mbox{ for any }X\in\GL_{2n}({\R_{\kaut^2}}),\label{eqn:Gtheta}\\
tr^{\sigma^{2}}t^{-1}&=tr^{\kaut^{-2}}t^{-1}=r &&\mbox{ for any }r\in{\R_{\kaut^2}}.\label{eqn:tsigma^2}
\end{align}
It then follows that we have $\theta^2=\id$.
Namely, for any $X\in \GL_{2n}({\R_{\kaut^2}})$,  
\begin{align}
X^{\theta^2} 
& = B^{-1}\ {}^t\!\left(B^{-1}\ {}^tX^{-\sigma^{-1}}\ B\right)^{-\sigma^{-1}}\ B\label{eqn:theta order 2}\\
& = B^{-1}\ {}^tB^{\sigma^{-1}}\  X^{\sigma^{-2}}\ {}^tB^{-\sigma^{-1}} B\nonumber\\
& = B^{-2} X^{\sigma^{-2}} B^{2}\nonumber\\
& = t^{-1}I_{2n} X^{\kaut^2} t I_{2n}\nonumber\\
& = X.\nonumber
\end{align}
We also have
\begin{align}\label{eqn:DetGtheta}
\Det(X^\theta)&=\Det(X)^{-\sigma^{-1}}.
\end{align}
Namely,  it is clear from \eqref{eq:adj mat c} and the fact that $\Det$ is a homomorphism, that for matrices $X, Y$ we have
 $\Det((XY)^\theta)=\Det(X^\theta Y^\theta)=\Det(X^\theta) \Det(Y^\theta)$.
Moreover, if $X$ is a transvection matrix, then so is $X^\theta$. Therefore we only have to check that~\eqref{eqn:DetGtheta} holds for diagonal matrices with $n-1$ trivial entries. However, this is clear.

We will now construct an amalgam $\amL^\delta$  inside $\SL_{2n}({\R_{\kaut^2}})$ that is isomorphic to the amalgam $\amG^\delta$. 
Consider the following matrix:
\begin{align}\label{eqn:\Cyc2nR}
\Cyc=\Cyc_{{\R_{\kaut^2}},2n}=\left(\begin{array}{@{}c|c@{}}
0 & I_{2n-1} \\
\hline
t & 0 \\
\end{array}\right).
\end{align}
We now define the automorphism $\Phi_{{\R_{\kaut^2}},2n}$ of $\SL_{2n}({\R_{\kaut^2}})$ given by $X\mapsto \Cyc^{-1}X\Cyc$.
Also define the map $i\colon \SL_2(\fk)\to \SL_{2n}({\R_{\kaut^2}})$ by
 $$\begin{array}{rl}
A &\mapsto \left( \begin{array}{c|c} 
A & \\
\hline
 & I_{2n-2}
 \end{array}\right)
 \end{array}.$$
Next, for $m=1,\ldots,n+1$, let 
$\phi_m\colon\SL_2(\fk)\to \SL_{2n}({\R_{\kaut^2}})$ by 
 $$\phi_m(A)=\Phi^{m-1}(i(A))\cdot \theta(\Phi^{m-1}(i(A)))$$
  and let $\L_m$ be the image of $\phi_{m}$.
  Note that 
\begin{align*}
\phi_{n+1}(A)=
\left( \begin{array}{cc|cc} {}^tA^{-\kaut^{-1}} &  \\
& I_{n-2}  \\ \hline & &  A & \\
 & & &  I_{n-2}\end{array} \right).
\end{align*}

Note that for each $m = 1,\dots n-1$  we have 
$$\L_{m}=\left\{\left( \begin{array}{ccc|ccc} I_{m-1} & & & & & \\
& A & & & & \\ & & I_{n-m-1} & & & \\ \hline & & &  I_{m-1} & &\\
 & & & & ^tA^{-\kaut} & \\  & & & & & I_{n-m-1}\end{array} \right) | A \in \SL_2(\fk)\right\}$$
and
$$\L_{n}=\left\{\left( \begin{array}{ccc|ccc} a^{\kaut^{-1}} & & & & & -t^{-1}b^\kaut\\
& I_{n-2} & & & & \\ & & a &b & & \\ \hline & &c &  d & &\\
 & & & & I_{n-2} & \\  -c^\kaut t & & & & & d^\kaut\end{array} \right) | \left(\begin{array}{cc}a& b\\ c& d\end{array}\right) \in \SL_2(\fk)\right\}.$$
 The latter can be verified more easily by observing that 
 \begin{align*}
{}^t\begin{pmatrix}a^{\kaut^{-1}}& -t^{-1}b^\kaut\\ -c^\kaut t& d^\kaut\end{pmatrix}^{-\sigma^{-1}}=
 \begin{pmatrix}a& -t^{-1}c^{\kaut^2}\\ -b^{\kaut^2}t& d^{\kaut^2}\end{pmatrix}^{-1}=
 \begin{pmatrix}d& t^{-1}c^{\kaut^2}\\ b^{\kaut^2}t& a^{\kaut^2}\end{pmatrix}. \end{align*}
One verifies that since $\Cyc^\theta=\Cyc$, we have $\theta\after\Phi=\Phi\after\theta$, and so for  $m=1,2,\ldots,n$, it follows that 
 \begin{align}\label{eqn:phi_k=Phi_{k-1}after phi_1}
 \phi_m=\Phi^{m-1}_{{\R_{\kaut^2}},2n}\after\phi_1.
 \end{align}
Let $I=\{1,2,\ldots,n\}$.
We shall denote the diagonal torus in the group $\L_{i}$ by $D_{i}$ for each $i\in I$.
For $(i,j)\ne (1,n)$ with $1\le i<j\le n$, let $\phi_{i,j}$ be the canonical isomorphism between 
 $\G_{i,j}=\langle \G_{i},\G_{j}\rangle$ and $\L_{i,j}=\langle \L_{i},\L_{j}\rangle_\G$ induced by $\phi_i$ and $\phi_{j}$. 
Moreover, let  $\phi_{n,1}$ be induced by $\phi_{n}$ and $\phi_{n+1}$.
It follows that $\L_{ij} \cong \SL_3(\fk)$ if $i-j\equiv  \pm {1 \mod  n}$ and $\G_{ij} \cong \L_{i}\times \L_{j} $ otherwise.

\bde
For each $i,j\in \{1,2,\ldots,n\}$, let 
 $\varphi_{i,j}\colon \L_{i}\into \L_{i,j}$ be the natural inclusion map.
Then we define the following amalgam:
 $$\amL^\delta=\{\L_{i},\L_{i,j},\varphi_{i,j}\mid i,j\in I\}.$$ 
 \ede
\bpr\label{prop:non-oriented CT amalgam}
The amalgam $\amL^\delta$ is contained in $\G^\theta$ and is isomorphic to $\amG^\delta$.
\epr
\bpf
That $\amL^\delta$ is contained in $\G^\theta$ follows by definition of $\phi_k$ and the fact that $\theta$ has order $2$ by~\eqref{eqn:theta order 2}.
We claim that the collection $\phi=\{\phi_i,\phi_{i,j}\mid i,j\in I\}$ is the required isomorphism between $\amG^\delta$ and $\amL^\delta$.
This is completely straightforward for all pairs $(i,j)$ except possibly for $(n,1)$. Here we have 
 $\phi_{n,1}\after \psi_{1,n}(A)=\phi_{n+1}(A^\delta)=\varphi_{1,n}\after \phi_1(A)$ since $A^\delta={}^t(A^\delta)^{-\kaut^{-1}}=A$.
\epf

\mn
We make some  observations on the form $\beta$ and the action of $\G$ on $\V$.
\ble\label{lem:trace valued form}
The form $\beta$ is non-degenerate trace-valued  and $(\sigma,t)$-sesquilinear. That is
 for all $u,v\in \V$ we have
 $\beta(v,u)=t\beta(u,v)^\sigma$
 and there exists $x\in {\Q_{\kaut^2}}$ such that
  $\beta(u,u)=x+tx^\sigma $.
\ele
\bpf
That $\beta$ is non-degenerate follows since $B$ is invertible.
To prove the second claim, let $u=\sum_{i=1}^n \lambda_i e_i+\mu_i f_i$ and
 let $u'=\sum_{i=1}^n \lambda_i' e_i+\mu_i' f_i$.
Using~\eqref{eqn:tsigma^2}, we find that 
\begin{alignat*}{2}
t\beta(u,u')^\sigma =t\left(\sum_{i=1}^n \lambda_i\mu_i'^\sigma  + \mu_it\lambda_i'^\sigma \right)^\sigma &=\sum_{i=1}^n t\mu_i'^{\sigma^2}\lambda_i^\sigma  + t\lambda_i'^{\sigma^2} t^{-1}\mu_i^\sigma\\
&=\sum_{i=1}^n \mu_i't\lambda_i^\sigma  + \lambda_i'\mu_i^\sigma =\beta(u',u).
                \end{alignat*}
Setting $u=u'$ and  $x=\sum_{i=1}^n \lambda_i\mu_i^\sigma$, and noting that $\mu_it\lambda_i^\sigma=t\mu_i^{\sigma^2}\lambda_i^\sigma=t(\lambda_i\mu_i^\sigma)^\sigma$, we get 
\begin{align*}
\beta(u,u)=\sum_{i=1}^n \lambda_i\mu_i^\sigma + \mu_it\lambda_i^\sigma
=x+tx^\sigma.
\end{align*}
\epf

\bde\label{dfn:dual basis adjoint map}
Given a ${\Q_{\kaut^2}}$-basis $\{a_1,\ldots,a_{2n}\}$ for $\V$, the {\dfn right dual basis } for $\V$ with respect to $\beta$ is the unique basis  $\{a_1^*,\ldots,a_{2n}^*\}$ such that $\beta(a_i,a_j^*)=\delta_{ij}$
 (note the order within $\beta$).
\ede

\ble\label{lem:dual standard basis}
If $\{a_1,\ldots,a_n,a_{n+1},\ldots,a_{2n}\}$ is a basis for $\V$ with Gram matrix $B$, then its right-dual basis is 
 $\{a_{n+1},\ldots,a_{2n},ta_1,\ldots,ta_n\}$.
\ele

\ble\label{lem:adjoint action}
If $g\in\GL(\V)$ is represented with respect to $\{a_1,\ldots,a_{2n}\}$ as right multiplication by a matrix $(g_{ij})$, then
 $g^*$ is represented with respect to $\{a_1^*,\ldots,a_{2n}^*\}$ as right multiplication by matrix ${}^t(g_{ij}^{\sigma^{-1}})^{-1}$. 
\ele
\bpf
Let $g^*$ be represented by $(g^*_{m,j})$. 
Then, 
\begin{align*}
\delta_{i,m}=\beta(a_i,a_m^*)&=\beta(g(a_i),g^*(a_m^*))\\
&=
\beta\left(\sum_j g_{i,j}a_j,\sum_j g^*_{m,j}a_j^*\right)
=\sum_j g_{i,j}(g^*_{m,j})^\sigma
\end{align*}
and so $(g_{i,j})\cdot {}^t(g_{j,m})^\sigma=I_{2n}$.
\epf

\bco\label{cor:duality preserves sM}
The right dual of an ${\R_{\kaut^2}}$-basis for $\M$ is an ${\R_{\kaut^2}}$-basis for $\M$.
\eco
\bpf
This follows from Lemmas~\ref{lem:dual standard basis} and~\ref{lem:adjoint action} by noting that $\GL(\M)$ is transitive on such bases and invariant under  $(g_{ij})\mapsto {}^t(g_{ij}^{\sigma^{-1}})^{-1}$.
\epf

\subsection{The geometry \texorpdfstring{$\De^\theta$}{} for \texorpdfstring{$\G^\theta$}{}}\label{subsection:flag transitive}
We now describe a geometry $\De^\theta$.
We shall subsequently prove that $\De^\theta$ is simply-connected, that $\G^\theta$ acts flag-transitively on $\De^\theta$, and that the amalgam of parabolic subgroups with respect to this action is the amalgam $\amB$ related to $\amL^\delta$ as in Proposition~\ref{prop:efficient completion recognition}.

Let $\De$ be the twin-building for the group $\G=\SL_{2n}({\R_{\kaut^2}})$ with twinning determined by $\M$ (for a construction see Subsection~\ref{subsec:building over R}).
Let $(W,S)$ be the Coxeter system with diagram $\tGamma$ of type $\tA_{2n-1}$. Call $S=\{s_i\mid i\in \tI\}$. 

\bde\label{dfn:theta on Delta}
For each $\cO_\vep$-lattice $\Lat_\vep$ we let 
 $$\Lat_\vep^\theta=\{v\in \V\mid \beta(u,v)\in \cO_\vep\mbox{ for all } u\in \Lat_\vep\}.$$
\ede

\ble\label{lem:tau on lattices}\ 
\begin{enumerate}
\item
If $\{a_1,\ldots,a_{2n}\}$ is a basis for $\V$ with right dual $\{a_1^*,\ldots,a_{2n}^*\}$ with respect to $\beta$, then
 $\Lat_\vep^\theta(a_1,\ldots,a_{2n})=\Lat_{-\vep}(a_1^*,\ldots,a_{2n}^*)$.
\item
For all $i$,  $j$ we have $(t^ja_i)^*=t^ja_i^*$ so 
  \begin{align*}
  \Lat_\vep^\theta(t^{j_1}a_1,\ldots,t^{j_{2n}}a_{2n})=\Lat_{-\vep}(t^{j_1}a_1^*,\ldots,t^{j_{2n}}a_{2n}^*).
  \end{align*}
 \item
 $\theta$ reverses inclusion of lattices.
\item
$\Lat_\vep^{\theta^2}(a_1,\ldots,a_{2n})=\Lat_\vep(ta_1,\ldots,ta_{2n})$.
\item
$\LatU_\vep^{\theta^2}(a_1,\ldots,a_{2n})=\LatU_\vep(a_1,\ldots,a_{2n}).$
 \end{enumerate}
\ele
\bpf
Parts (a) and (b) are straightforward consequences of the fact that $\beta$ is $\sigma$-sesquilinear.
Part (c)  follows from Definition~\ref{dfn:theta on Delta}. 
Part (d) and (e): By Lemma~\ref{lem:trace valued form}, we have $\beta(u,v)=t\beta(v,u)^\sigma\in \cO_\vep$, so the right dual basis of $\{a_1^*,\ldots,a_{2n}^*\}$ is $\{ta_1,\ldots,ta_{2n}\}$ and the claim follows from (a).
\epf

\mn
The standard chamber in $\De_\vep$ is  $c_\vep(e_1,\ldots,e_n,f_1,\ldots,f_n)$.

\bpr\label{prop:tau on Delta}
The map $\theta$ is an involution on $\De$ that induces isomorphisms
 $\theta\colon \De_\vep\to \De_{-\vep}$ where $\typ(\theta)\colon \tI\to \tI$ is the graph
  isomorphism defined by $i\to i-n  \mod (2n)$.
Moreover, $\theta$ interchanges the standard chambers $c_+$ and $c_-$.
\epr
\bpf
By Lemma~\ref{lem:tau on lattices} (a) and (c) $\theta$ sends admissible chains of $\cO_\vep$-lattices to admissible chains of $\cO_{-\vep}$-lattices.
In particular, it interchanges $\De_\vep$-objects with $\De_{-\vep}$-objects while preserving incidence. 
Thus $\theta$ induces the required isomorphisms. 
By  Lemma~\ref{lem:tau on lattices} (d) $\theta$ is an involution.
 We now analyze how types are permuted by $\theta$.

Let $C_{i,\vep}$ be the object of type $i$ on $c_\vep$.
We show that $C_{i,\vep}^\theta=C_{n+i,-\vep}$.
This follows immediately from Lemmas~\ref{lem:tau on lattices}~and~\ref{lem:dual standard basis}.
In particular $c_+$ and $c_-$ are interchanged.

Let $d_\vep\in \De_\vep$ be any other chamber. Then, since $\SL_{2n}({\R_{\kaut^2}})$ is transitive on chambers of $\De_\vep$, it contains an element $g$ such that 
 $g(c_\vep)=d_\vep$.
By Corollary~\ref{cor:duality preserves sM}, ${}^t(g^{\sigma^{-1}})^{-1}$ takes $c_{-\vep}$ to a chamber
 $d_{-\vep}$ that is opposite to $d_\vep$ and such that 
  $(gd_\vep)^\theta=d_{-\vep}$.
As $v_\vep(\Det(g))=v_\vep\Det(({}^t(g^{\sigma^{-1}})^{-1}))$, (where $\Det$ denotes the Dieudonn\'e  determinant), $\theta$ permutes the types on $d_\vep$ as it does on $c_\vep$.
\epf

\bde\label{dfn:theta on W}
We shall abuse notation and write $\theta(i)=\typ(\theta)(i)=i-n$ for $i\in \tI$.
Thus $\theta$ is a graph automorphism of $\tGamma$ inducing an automorphism of the Coxeter system $(W,S)$, which we shall also denote $\theta$.
\ede

\bde\label{dfn:Detau}
We define a relaxed incidence relation on $\De_\vep$ as follows.
We say that $d_\vep$ and $e_\vep$ are $(i,\theta(i))$-adjacent if and only if 
$d_\vep$ and $e_\vep$ are in a common $\{i,\theta(i)\}$-residue.
In this case we write 
$$d_\vep\approx_i e_\vep,$$
where we let $i\in I=\{1,\ldots,n\}$.
Note that the residues in this chamber system are $J$-residues of $\De_\vep$ where $J^\theta=J$.
In Subsection~\ref{sec:simple-connectedness}
 we shall see that the resulting chamber system $(\De_\vep,\approx)$ is simply connected.
Let
 $$\De^\theta=\{(d_+,d_+^\theta)\mid d_+\opp d_+^\theta\}.$$
Adjacency is given by $\approx$. It is easy to see that residues of $\De^\theta$ are the intersections of residues of $(\De,\approx)$ with the set $\De^\theta$.
\ede

\ble\label{lem:De(tau) by bases}
$(d_+,d_-)\in\De^\theta$ if and only  if there exists  $\{a_1,\ldots,a_n,b_1,\ldots,b_n\}$, an ${\R_{\kaut^2}}$-basis for $M$ whose Gram matrix is $B$ and  $d_\vep = c_\vep(a_1,\ldots,a_n,b_1,\ldots,b_n)$ for $\vep=+,-$.
\ele
\bpf
As in the proof of Proposition~\ref{prop:tau on Delta}, one verifies that any such basis gives rise to a pair of chambers in $\De^\theta$.
Conversely, let $(d_+,d_-)\in \De^\theta$.
That means that $d_-=d_+^\theta$.
Let $\Sigma=\Sigma(d_+,d_-)$ be the twin-apartment containing $d_+$ and $d_-$.
Then $\Sigma^\theta=\Sigma$.
Let $\{a_1,\ldots,a_n,b_1,\ldots,b_n\}$ be an ${\R_{\kaut^2}}$-basis for $\M$ such that
 $\Sigma=\Sigma\{a_1,\ldots,a_n,b_1,\ldots,b_n\}$ and
 $d_\vep=c_\vep(a_1,\ldots,a_n,b_1,\ldots,b_n)$, where
  $\langle a_1,\ldots,a_n,b_1,\ldots,b_n\rangle_{\cO_\vep}$ has type $0$.
Let $\{a_1^*,\ldots,a_n^*,b_1^*,\ldots,b_n^*\}$ be the right dual basis with respect to $\beta$.
Then, since $\{d_+^\theta,d_-^\theta\}=\{d_+,d_-\}$ uniquely determines $\Sigma$, it follows from Lemma~\ref{lem:tau on lattices}, that, for $\vep=\pm$,  
$$\begin{array}{ll}
\Sigma&=\Sigma\{a_1^*,\ldots,a_n^*,b_1^*,\ldots,b_n^*\},\\
d_\vep&=c_\vep(a_1^*,\ldots,a_n^*,b_1^*,\ldots,b_n^*).\\
\end{array}$$
By Corollary~\ref{cor:duality preserves sM} both bases are ${\R_{\kaut^2}}$-bases for $\M$.
Note that the type of the lattice 
$\langle a_1^*,\ldots,a_n^*,b_1^*,\ldots,b_n^*\rangle_{\cO_\vep}=\langle a_1,\ldots,a_n,b_1,\ldots,b_n\rangle_{\cO_{-\vep}}^\theta$ is $n$.
Now consider the ${\R_{\kaut^2}}$-linear map
 $$\begin{array}{rl}
 \phi\colon \M &\to \M\\
  b_i & \mapsto a_i^*\\
  ta_i & \mapsto b_i^*\\
 \end{array}$$
for all $i=1,2,\ldots,n$.
It is easy to check that $\phi$ is a type-preserving automorphism of $\De_\vep$ such that
 $d_\vep^\phi=d_\vep$ since it is an ${\R_{\kaut^2}}$-linear map that sends the object of type $i$ on $d_\vep$
  to the object of type $i$ on $d_\vep$.
It follows from Lemma~\ref{lem:stabilizer of c+ and c-} that
  $$\begin{array}{ll}
  b_i & =\lambda_i t^ka_i^*,\\
  ta_i &= \mu_i t^kb_i^*,\\
 \end{array}$$
where $\lambda_i,\mu_i\in \fk^*$ and $k\in \ZZ$.
Computing $\beta(b_i,b_i^*)$ and using that $\beta(a^*_i,ta_i)=1$, we find $k=0$ and $\mu_i=\lambda_i^{\sigma^{-1}}$.
Without modifying the chambers $d_\vep$, we may replace $a_i$ by $\lambda_i^{-\sigma}a_i$ and keep $b_i$ so that 
  $$\begin{array}{ll}
  b_i & =a_i^*,\\
  ta_i &=b_i^*,\\
 \end{array}$$
and so the Gram matrix of $\{a_1,\ldots,a_n,b_1,\ldots,b_n\}$ is $B$.
\epf

\mn
Let $\GUD=\GD\cap \GU_n({\R_{\kaut^2}})$ and $\SUD=\GD\cap \SU_n({\R_{\kaut^2}})$.
\ble\label{lem:GUD}\ 
\begin{enumerate} 
\item  $\GUD=\{\diag(\lambda_1,\ldots,\lambda_n,\lambda_1^{-\sigma^{-1}}, \ldots, \lambda_n^{-\sigma^{-1}})\mid \lambda_1,\ldots,\lambda_n\in {\R_{\kaut^2}^*}\}$,

\item If $\fk/\fk_\kaut$ is a cyclic Galois extension, then 
$\Det$ is onto and $N_{{\R_{\kaut^2}^*}/\A^*}$ is an isomorphism:
\begin{align*}
\GUD&\stackrel{\Det}{\longrightarrow} 
\{at^m[{\R_{\kaut^2}^*},{\R_{\kaut^2}^*}]\mid a\in \ker N_{\fk/\fk_\kaut}, m\in \ZZ\}\\
&\stackrel{N_{{\R_{\kaut^2}^*}/\A^*}}{\longrightarrow}
\{b((-1)^{s/2-1}T)^m\mid b\in \n_{\kaut^{\scriptscriptstyle 2}}\cap \ker N_{\fk_{\kaut^2}/\fk_{\kaut}}, m\in \ZZ\}.
\end{align*}

\item If $N_{\fk_{\alpha^2}/\fk_\alpha}$ is surjective, then, 
$\SUD=\Dt=\langle \phi_i(D_{i})\mid i\in I\rangle$.\\
Moreover, 
\begin{align*}\Dt=\{\diag(\lambda_1,\ldots,\lambda_n,\lambda_1^{-\sigma^{-1}}, \ldots,\lambda_n^{-\sigma^{-1}})\mid &\lambda_1,\ldots,\lambda_n\in \fk^*, \\
&
\prod_{i=1}^n \lambda_i\lambda_i^{-\kaut}\in \ker N_{\fk/\fk_{\kaut^2}}\}.
\end{align*}
\end{enumerate}
\ele
\bpf
(a) Let $\psi\in \GD$. By Lemma~\ref{lem:stabilizer of c+ and c-} this means that 
$$\begin{array}{rl}
 \psi\colon \M &\to \M\\
  e_i & \mapsto \lambda_i t^m e_i\\
  f_i & \mapsto \mu_i t^m f_i\\
 \end{array}$$
with $\lambda_i,\mu_i\in \fk$ for all $i=1,2,\ldots,n$ and some $m\in \ZZ$.

The conditions $\beta(\lambda_i e_i,\mu_j f_j)=\delta_{ij}$ (and, equivalently $\beta(\mu_j f_j,\lambda_i e_i)=t\delta_{ji}$) yield $\mu_i=\lambda_i^{-\kaut}$, but no restriction on $k$. Any such element lies in $\GUD$.

(b) From (a) we find that $\Det(\psi)=b=c^{\kaut}c^{-1}$, where $c=\prod_{i=1}^n \lambda_i^{-1}$. Clearly any $b$ of this form appears as $\Det(\psi)$ of some $\psi$.
By Hilbert's 90th theorem, therefore $\Det$ is onto.

Note that by Lemma~\ref{lem:norms}, the map $N_{{\R_{\kaut^2}^*}/\A^*}$ is injective. It suffices therefore to check that this restriction is onto. 
First note that it sends $t\mapsto (-1)^{s/2-1} T$. 
To check that its restriction $N_{\fk/\fk_{\kaut^2}}$ is onto, consider the  following diagram: 
\begin{align*}
\xymatrix{
\fk^*\ar@{->>}[rr]_{N_{\fk/\fk_{\kaut^2}}} \ar@(ur,ul)@{->>}[rrrr]^{N_{\fk/\fk_\kaut}} && \n_{\kaut^2}\ar@{->>}[rr]_{N_{\fk_{\kaut^2}/\fk_{\kaut}}} && \n_\kaut.
}
\end{align*}

Note that all maps are surjective since $N_{\fk/\fk_\kaut}=N_{\fk_{\kaut^2}/\fk_{\kaut}}\after N_{\fk/\fk_{\kaut^2}}$.
It follows that $N_{\fk/\fk_{\kaut^2}}$ takes 
 $\ker N_{\fk/\fk_\kaut}$ to $\n_{\kaut^2}\cap \ker N_{\fk_{\kaut^2}/\fk_\kaut}$.
 
(c) It is clear from the definition of the $\L_{i}$ that $\Dt\le \SUD$.
With $\psi$ as in (a) we find that $m=0$ and 
 $\Det(\psi)\in [{\R_{\kaut^2}^*},{\R_{\kaut^2}^*}]=\ker N_{\fk/\fk_{\kaut^2}}$ by Lemma~\ref{lem:norms} and Hilbert's 90th theorem.

To see $\SUD\le \Dt$, let  
 $\psi=\diag(\lambda_1,\ldots,\lambda_n,\lambda_1^{-\sigma^{-1}}, \ldots,\lambda_n^{-\sigma^{-1}}) \in \SUD$, that is, $\lambda_1,\ldots,\lambda_n\in \fk^*$ and  $\prod_{i=1}^n \lambda_i\lambda_i^{-\kaut}=d^{-1}d^{\kaut^2}\in \ker N_{\fk/\fk_{\kaut^2}}$.
Let 
\begin{align*}
\eta=\phi_n\left(\begin{pmatrix} d^{-\kaut} & 0 \\ 0 & d^\kaut\end{pmatrix}\right).
\end{align*}
Then, $\eta^{-1}\psi$ is a diagonal matrix of determinant $1$.
Now suppose that 
$\eta^{-1}\psi=\diag(\mu_1,\ldots,\mu_n,\mu_1^{-\kaut},\ldots,\mu_n^{-\kaut})$ such that 
 $\prod_{i=1}^n \mu_i\mu_i^{-\alpha}=1$.
Let $a=\prod_{i=1}^n \mu_i$. Then, $a=a^{\kaut}$, so $a\in \fk_{\kaut}$. By assumption there exists some $c\in \fk_{\kaut^2}$ with $cc^\kaut=a$.
Let 
\begin{align*}
\gamma=\phi_n\left(\begin{pmatrix} c^{-\kaut} & 0 \\ 0 & c^\kaut\end{pmatrix}\right).
\end{align*}
Then, $\gamma\eta^{-1}\psi\in \langle \phi_i(D_{i})\mid i\in\{1,2,\ldots,n-1\}\rangle$.
This shows that $\psi\in \Dt$.
\epf

\bth\label{thm:flag-transitive}
Assume that $\fk/\fk_\kaut$ is cyclic and Galois. The group $\G^\theta$ acts flag-transitively on $\De^\theta$.
\eth
\bpf
Let $(d_+,d_-)\in \De^\theta$.
By Lemma~\ref{lem:De(tau) by bases} there exists
$\{a_1,\ldots,a_n,b_1,\ldots,b_n\}$,  an ${\R_{\kaut^2}}$-basis for $\M$  with Gram matrix $B$.
The ${\R_{\kaut^2}}$-linear map
 $$\begin{array}{rl}
 x\colon \M &\to \M\\
  e_i & \mapsto a_i\\
  f_i & \mapsto b_i\\
 \end{array}$$
for all $i=1,2,\ldots,n$
belongs to $\GU_{\R_{\kaut^2}}(\M,\beta)$ and sends $(c_+,c_-)$ to $(d_+,d_-)$.
Now suppose $x$ is represented by $X\in \GL_{2n}({\R_{\kaut^2}})$ and let $a$ represent $\Det(G)$ in ${\R_{\kaut^2}^*}/[{\R_{\kaut^2}^*},{\R_{\kaut^2}^*}]$.
As $X$ preserves types, $v_\vep (\Det(G))=2nm$ for some $m\in \ZZ$ and since $(t^{-m} X)^\theta=t^{-m} X^\theta$ we may assume $v_\vep(\Det(X))=0$, so that $a\in \fk$.
Then, by~\eqref{eqn:DetGtheta} we have 
\begin{align*}
aa^{\sigma^{-1}}=aa^\kaut\in [{\R_{\kaut^2}^*},{\R_{\kaut^2}^*}].
\end{align*}
By Lemma~\ref{lem:norms}, $aa^\kaut=c^{\kaut^2}c^{-1}$ for some $c\in \fk$.
Hence 
\begin{align*}
N_{\fk/\fk_\kaut}(a)=N_{\fk/\fk_{\kaut^2}}(aa^\kaut)=N_{\fk/\fk_{\kaut^2}}(c^{\kaut^2}c^{-1})=1.
\end{align*}
By Lemma~\ref{lem:GUD} there is $y\in \GUD$ such that $y\after x\in \SU_{\R_{\kaut^2}}(\M,\beta)$. Clearly also $y\after x$ takes $(c_+,c_-)$ to $(d_+,d_-)$, as desired.
\epf

\subsection{Simple connectedness}\label{sec:simple-connectedness}
In this subsection we will prove that 
the chamber system $(\De^\theta,\approx)$ is connected and simply-connected.
In order to do so we shall in fact prove a stronger result, namely that
 $(\De^\theta,\sim)$ is connected and simply connected.
Namely,
\ble\label{lem:sim to approx}
Suppose that $X$ is a subset of $\De_+$ such that
 $(X,\sim)$ is connected and simply connected. Then $(X,\approx)$ is also connected and simply connected.
\ele
\bpf
Note that each rank $r<n$ residue of $(\De_+,\sim)$  is included in a residue of rank $\le r$ of $(\De_+,\approx)$.
Since connectedness is a statement about rank $1$ residues and simple connectedness is a statement about rank $2$ residues, we are done.
\epf

\mn We will use the techniques developed in~\cite{DevMuh2007} to show that $(\De^\theta,\sim)$ is simply connected.
\bde\label{dfn:residual filtration} In the terminology of loc. cit.\ a collection $\{\cC_m\}_{m\in \NN}$ of subsets of a chamber system $\cD$ over $I$ is a {\dfn filtration} if
 the following are satisfied:
\begin{itemize}
\item[F1] For any $m\in \NN$ $\cC_m\sbe \cC_{m+1}$,
\item[F2] $\bigcup_{m\in \NN} \cC_m=\cD$,
\item[F3] For any $m\in \NN_{>0}$, if $\cC_{m-1}\ne \emptyset$, there exists an $i\in I$ such that for any
 $c\in \cC_m$, there is a $d\in \cC_{m-1}$ that is $i$-adjacent to $c$.
\end{itemize}
It is called a {\dfn residual filtration} if the intersections of $\cC$ with any given residue is a filtration of that residue.
\ede

For any $c\in \cD$, let $|c|=\min\{\lambda\mid c\in \cC_\lambda\}$. For a subset $X\sbe \cD$ we accordingly define
\begin{align*}
 |X|=\min\{|c| \mid c\in X\} \mbox{ and } \\\aff(X)=\{c\in X\mid |c|=|X|\}.
\end{align*}
We shall make use of the following result from~\cite{DevMuh2007}.
\bth\label{thm:DeMu}~\cite[Theorem 3.14]{DevMuh2007}
Suppose $\cC$ is a residual filtration on $\cD$ such that for any rank $2$ residue $R$, $\aff(R)$ is connected and for any rank $3$ residue $R$, $\aff(R)$ is simply $2$-connected, then the following are equivalent.
\begin{enumerate}
\item $\cD$ is simply $2$-connected.
\item $\cC_n$ is simply $2$-connected for all $n\in \NN$.
\end{enumerate}
\eth
We now let $\cD$ be the chamber system $\De_+$, with adjacency relations $\approx_i$ ($i\in I$).
We then define a residual filtration $\cC$ on $\De_+$ with the property that 
 $\cC_0\cong\De^\theta$. We shall use that $\De_+$ is simply connected. In order to obtain simple connectedness of $\De^\theta$ it  will suffice to show that $\cC$ satisfies the conditions of the theorem.

\subsection{The filtration \texorpdfstring{$\cC$}{}}\label{subsec:filtration}
Recall that $(W,S)$ is a Coxeter system with diagram $\tGamma$ of type $\tA_{2n-1}$, where $S=\{s_i\mid i\in \tI\}$. 
For any $w\in W$, let $l(w)$ denote its length with respect to $S$. Recall from Definition~\ref{dfn:theta on W} that $\theta$ acts on $\tI$ and $(W,S)$.
In order to define the filtration $\cC$ we first let 
 $$\delta^\theta(W)=\{w\in W\mid \exists d_\vep\in \De_\vep\colon w=\delta_*(d_\vep,d_\vep^\theta)\}.$$
We also fix an injective map $|\cdot|\colon \delta^\theta(W)\to \NN$ such that whenever
 $l(w)>l(w')$, we have $|w|>|w'|$ and $|1|=0$. For any $m\in \NN$, 
we then define a filtration on $\De_+$ using $|\cdot|$ as follows:
Let   $$\cC_m=\{c_+\in\De_+\mid |\delta_*(c_+,c_+^\theta)|\le m\}.$$
In particular we have 
 \begin{align}\label{eqn:C_0}
\cC_0=\{c_+\in \De_+\mid (c_+,c_+^\theta)\in \De^\theta\}.
 \end{align}
In fact the map $(\De^\theta,\sim)\to (\cC_0,\sim)$ sending $(d_+,d_+^\theta)\mapsto d_+$ is an isomorphism of chamber systems.

In the remainder of this section we prove that $\cC$ is a residual filtration. First however, we will need some technical lemmas about $\delta^\theta(W)$.
Let
 \begin{align*}
 \Inv^\theta(W)&=\{u\in W\mid u^\theta=u^{-1}\},\\
W(\theta)&=\{w(w^{-1})^\theta \mid w\in W\}.
\end{align*}

 These elements are called twisted involutions in~\cite{Spr1985} and
 ~\cite{GraHorMuh2011}. Some of the results below have somewhat weaker
 forms in the most general case of a quasi-twist. See
 ~\cite{GraHorMuh2011} for details on both twisted involutions and of
 the corresponding geometries.

We now have the following:
\ble\label{lem:Wtau}
 $$ \Inv^\theta(W)=W(\theta).$$
More precisely, given any $u\in \Inv^\theta(W)$ there exists a word $w\in W$ such that
 $w(w^{-1})^\theta$ is a reduced expression for $u$.
\ele
\bpf
Clearly we have $W(\theta)\sbe \Inv^\theta(W)$.
Let $w\in \Inv^\theta(W)$. 
Then, by~\cite[Proposition 4.3]{GraHorMuh2011} or ~\cite[Proposition 3.3(a)]{Spr1985} there exists a spherical subset $J\sbe \tI$ and $s_1,\ldots,s_h\in S$ such that 
 $w=s_1\cdots s_h w_J s_{\theta{h}}\cdots s_{\theta{1}}$, where $w_J$ denotes the longest word in $W_J$.
Note since $J$ is spherical and simply-laced, $\tGamma_J$ has a $\theta$-fixed vertex or edge or $\tGamma_J=\tGamma_{J_1}\uplus \tGamma_{J_1^\theta}$, for some $J_1\sbne J$.
Since $\theta$ has no fixed points or edges on $\tGamma$, we are in the latter case.
Hence $w_J=w_{J_1}\cdot w_{J_1^\theta}=
 w_{J_1}\cdot w_{J_1^\theta}\in  W(\theta)$ and so $w\in W(\theta)$.
 \epf

\bre
Note that the proof of Lemma~\ref{lem:Wtau} only uses that the diagram is simply-laced and the involution $\theta$ has no fixed nodes or edges.
\ere

\ble\label{lem:tau does not preserve roots}
$\theta$ does not commute with any reflection.
\ele
\bpf
Let $r$ be any reflection such that $r^\theta=r$.
Then in fact $r\in \Inv^\theta(W)=W(\theta)$. However, all elements of $W(\theta)$ have even length and $r$ being a conjugate of a fundamental reflection does not.
\epf
\ble\label{lem:sws'<>w}
For $u\in \Inv^\theta(W)$ and $i\in \tI$, we have
 $l(s_ius_{\theta(i)})= l(u)\pm 2$.
\ele
\bpf
Suppose that $l(s_ius_{\theta(i)})=l(u)$, then by Lemma 4.2 of ~\cite{GraHorMuh2011} $s_ius_{\theta(i)}=u$,  contradicting Lemma~\ref{lem:tau does not preserve roots}.
\epf
%

\mn The following lemma characterizes $\delta^\theta(W)$.
\ble\label{lem:Wtau=d_*}
 $\delta^\theta(W)=\Inv^\theta(W)$.
\ele
\bpf
Let $c_\vep\in \De_\vep$.
Then $u=\delta_*(c_\vep,c_\vep^\theta)$ satisfies $u^\theta=u^{-1}$. Therefore the inclusion
 $\sbe$ follows  by definition.
Conversely, consider a chamber $c_\vep$ such that $c_\vep\opp c_\vep^\theta$.
Then the apartment $\Sigma(c_\vep,c_\vep^\theta)$ is preserved by $\theta$ and identifying it with the
Coxeter group we see that $\theta$ acts on $\Sigma$ as it acts on $W$.
Let $u\in \Inv^\theta(W)$. Then, by Lemma~\ref{lem:Wtau} it is of the form $w (w^{-1})^\theta$ for some $w\in W$.
Let $d_\vep$ be the chamber such that $\delta_\vep(c_\vep,d_\vep)=w$, then by induction on the length $l(w)$ and Lemma~\ref{lem:sws'<>w} we have 
 $\delta_*(d_\vep,d_\vep^\theta)=w(w^{-1})^\theta=u$ as desired.
\epf

\mn
In the sequel we shall use the following notation for projections. Given a residue $R$ of $\De_\vep$, we denote projection from $\De_\vep$ onto $R$ by $\proj_R$ and denote (co-) projection from $\De_{-\vep}$ onto $R$ by $\proj_R^*$.

\ble\label{lem:equal distances}
Suppose that $c_\vep\in \De_\vep$ satisfies $\delta_*(c_\vep,c_\vep^\theta)=w$, let $i\in \tI$ and suppose that
 $\pi$ is the $\sim_i$-panel on $c_\vep$.
 Then,
 \begin{enumerate}
 \item If $l(s_iw)>l(w)$, then all chambers $d_\vep\in \pi-\{c_\vep\}$ except one satisfy
 $\delta_*(d_\vep,d_\vep^\theta)=w$.
The last chamber $\check{c}_\vep$ satisfies $\delta_*(\check{c}_\vep,(\check{c}_\vep)^\theta)=s_iws_{\theta(i)}$.
 \item If $l(s_i w)<l(w)$, then all chambers  $d_\vep\in \pi-\{c_\vep\}$ have the propery
 $\delta_*(d_\vep,d_\vep^\theta)=s_i w s_{\theta(i)}$.
 \end{enumerate}
In particular, if $w=1$, then all chambers $d_\vep\in \pi-\{c_\vep\}$ except one satisfy
 $\delta_*(d_\vep,d_\vep^\theta)=1$.
\ele

\bpf
This follows from Lemma 4.6~\cite{GraHorMuh2011} and Lemma~\ref{lem:sws'<>w}.
\epf

\mn
We define the following subset of a given $J$-residue $R$:\begin{align}\label{eqn:A_theta(R)}
A_\theta(R)=\{c\in R\mid l(\delta_*(c,c^\theta))\mbox{ is minimal among all such distances}\}.
\end{align}
In particular, if $R\opp_\De R^\theta$, then 
 \begin{align}\label{eqn:A_theta(R) in De^theta}
 A_\theta(R)=\{c\in R\mid (c,c^\theta)\in \De^\theta\}.
 \end{align}
\ble\label{lem:A_tau has unique distance}
Let $R$ be a $J$-residue of $\De_\vep$.
Let $c\in A_\theta(R)$, $w=\delta_*(c,c^\theta)$ and let $d \in R$.
Then, $d\in  A_\theta(R)$ if and only if $w=\delta_*(d,d^\theta)$.
Moreover, $w$ is determined by the fact that for any $j\in J$ we have
 $l(s_jw)=l(w)+1$.
\ele
\bpf
First note that by Lemma~\ref{lem:equal distances}, $\{\delta_*(x,x^\theta)\mid x\in R\}=\{uwu^\theta\mid u\in W_J\}$.
Moreover, the coset $W_JwW_{\theta(J)}$ has a minimal element $m$ that is characterized by the fact that
 $l(s_jm)=l(m)+1$ and $l(ms_{\theta(j)})=l(m)+1$ for all $j\in J$.
We claim that $w$ has that property as well.
Namely, let $j\in J$ have the property that $l(ws_{\theta(j)})=l(s_jw)<l(w)$.
Then, by Lemma~\ref{lem:equal distances} (b) any element $d$ in the $j$-panel on $c$ has the property that $\delta_*(d,d^\theta)=s_jws_{\theta(j)}$ and by Lemma~\ref{lem:sws'<>w} this must have length $l(w)-2$, a contradiction to the fact that $c\in A_\theta(R)$.
Thus, $w$ satisfies the conditions on $m$ and it follows that $w=m$.
\epf

\bpr\label{prop:A_tau has unique distance}
Let $c\in R$ and let $w=\delta_*(c,c^\theta)$.
The following are equivalent:
\begin{enumerate}
\item $c\in A_\theta(R)$.
\item $w=w_R$, the unique element of minimal length in $W_JwW_{\theta(J)}$.
\item $c\in \cC_k$, where $k=\min\{l\mid \cC_l\cap R\ne \emptyset\}$.
\end{enumerate}
In particular, we have $A_\theta(R)=\aff(R)$.
\epr
\bpf
By Lemma~\ref{lem:A_tau has unique distance} (a) and (b) are equivalent.
Since $|\cdot|$ is strictly increasing, also (b) and (c) are equivalent.
\epf

\bpr\label{prop:is filtration}
$\cC$ is a residual filtration.
\epr
\bpf
We check the conditions in Definition~\ref{dfn:residual filtration}.
Part (F1) and (F2) are immediate.
Now let $R$ be a $J$-residue and suppose that $R\cap\cC_{n-1}\ne \emptyset$.
If $R\cap \cC_n=R\cap \cC_{n-1}$ there is nothing to check, so assume otherwise and let $w\in \delta^\theta(W)$ be unique with $|w|=n$.
Also pick any $c\in R\cap \cC_n-\cC_{n-1}$ so that $w=\delta_*(c,c^\theta)$.
By Proposition~\ref{prop:A_tau has unique distance}, $c\not\in A_\theta(R)$ and so, by Lemma~\ref{lem:A_tau has unique distance},
 there exists a $j\in J$ with $l(s_jw)<l(w)$. 
Therefore by Lemma~\ref{lem:sws'<>w}, any $j$-neighbor $d$ of of $c$ has $l(\delta(d,d^\theta))=l(w)-2$ and therefore belongs to $\cC_{n-1}$.
\epf

\subsection{Simple connectedness of $\De^\theta$}
Proposition~\ref{prop:is filtration} allows us to apply Theorem~\ref{thm:DeMu} and, by Proposition~\ref{prop:A_tau has unique distance}, in order to show simple connectedness of $\De^\theta$, it suffices to show that $\aff(R)=A_\theta(R)$ is connected when $R$ has rank $2$ and is simply connected when $R$ has rank $3$.
We shall first obtain some general properties of $A_\theta(R)$ and then verify the connectedness properties using concrete models of $A_\theta(R)$.

\bpr\label{prop:AtauR}{\rm (See Corollary 7.4 of~\cite{BlHo2009})}
For $\vep=\pm$, let $S_\vep\sbne R_\vep$ be residues of $\De_\vep$ such that
 $S_\vep=\proj^*_{R_\vep}(R_{-\vep})$ and let $x_\vep\in R_\vep$ be an arbitrary chamber
 and assume in addition that
 $R_{-\vep}=R_\vep^\theta$ and $x_{-\vep}=x_\vep^\theta$, for $\vep=\pm$.
Then,
 $x_\vep\in A_\theta(R_\vep)$ if and only if
 \begin{enumerate}
 \item $x_\vep$ belongs to a residue opposite to $S_\vep$ in $R_\vep$ whose type is also opposite to the type of $S_\vep$ in $R_\vep$ and

 \item $\proj_{S_\vep}(x_\vep)\in A_\theta(S_\vep)$.
 \end{enumerate}
\epr
\bpf
This is exactly the same as the proof in~\cite{BlHo2009}~noting that 
 it suffices for $\theta$ to be an isomorphism between $\De_+$ and $\De_-$ that preserves lengths of codistances.
\epf

\mn
Recall that for a spherical residue $X_\vep \sbe \De_\vep$ and $x_{\vep},z_\vep \in \De_{\vep}$, the chamber $y_\vep=\proj^*_{X_\vep}(x_{-\vep})$ is the unique chamber in $X_\vep$ having maximal length codistance to $x_{-\vep}$. For all $z_\vep\in X_\vep$ it satisfies
\begin{align}\label{eqn:codistance}
\delta_*(z_\vep,x_{-\vep})= \delta_\vep(z_\vep,y_\vep)\delta_*(y_\vep,x_{-\vep}).
\end{align}

\ble\label{lem:inverse projections}
With the notation of Proposition~\ref{prop:AtauR}, 
 $\proj^*_{S_\vep}$, $\proj^*_{S_{-\vep}}$  define
 adjacency preserving bijections between $S_{-\vep}$ and $S_\vep$ such that $(\proj^*_{S_\vep})^{-1}=\proj^*_{S_{-\vep}}$.
Let $l=\max\{l(\delta_*(c_\vep,d_{-\vep}))\mid c_\vep\in S_\vep,d_{-\vep}\in S_{-\vep}\}$.
Then, 
$d_{-\vep}=\proj^*_{S_{-\vep}}(c_\vep)$ if and only if 
 $l(\delta_*(c_\vep,d_{-\vep}))=l$.
\ele
\bpf
This is the twin-building version of the main result of~\cite{DreSch1987}.
\epf

\mn
In view of Proposition~\ref{prop:AtauR}, in order to study $A_\theta(R)$ entirely inside $R$ we need to know what $A_\theta(S)$ looks like if $\proj^*_{S}\after\theta$ is a bijection on $S$.
From now on we shall write $\theta_S=\proj^*_{S}\after\theta$.
\bco\label{cor:proj after tau has order 2}
In the notation of Proposition~\ref{prop:AtauR}, 
$\theta_{S_\vep}$ has order $2$.
\eco
\bpf 
Let $c\in S_\vep$. 
Then 
  $l(\delta_*(c^\theta,(\proj^*_{S_{-\vep}}(c))^\theta))=l(\delta_*(c,(\proj^*_{S_{-\vep}}(c))))$.
Therefore, by Lemma~\ref{lem:inverse projections}, $\proj^*_{S_{\vep}}(c^\theta)=(\proj^*_{S_{-\vep}}(c))^\theta$.
The claim of the lemma follows.
\epf

\mn The next proposition describes the structure of the residues of $\De^\theta$.
\bpr\label{prop:theta residues}
Let $J\sbne \tI$ be $\theta$-invariant and suppose $R$ is a $J$-residue of $\De_+$ such that $(R,R^\theta)$ meets $\De^\theta$ in a residue of $\De^\theta$.
Then, 
\begin{enumerate}
\item\label{types}
 $J=J_1\uplus J_1^\theta$ and $\tGamma_J=\tGamma_{J_1}\uplus \tGamma_{{J_1^\theta}}$, 
\item\label{PQ}  $R=P\times Q^\theta$ and $R^\theta=P^\theta\times Q$, where $P\sbe R$ and $Q\sbe R^\theta$ are arbitrary ${J_1}$-residues,
\item\label{proj}
we can pick $P$ and $Q$ so that 
$\proj^*_R\colon P^\theta\to Q^\theta$ and $\proj^*_R\colon Q\to P$ are (possibly type changing) isomorphisms,
\item\label{Q=P^theta_R}
$R\cong P\times P^{\theta_R}$, where $P$ is a residue of type ${J_1}$,
\item\label{AthetaR} we have 
 $A_\theta(R)=\{(p,q)\in P\times P^{\theta_R}\mid p\opp_P q^{\theta_R}\}$. 
 In particular, $A_\theta(R)$ is isomorphic to the geometry of pairs of opposite chambers in $P$.

\end{enumerate}
\epr
\bpf
\eqref{types} Since $J\ne \tI$, there is $i\in \tI$ with $J\sbe \tI-\{i\}$, hence in fact $J\sbe \tI-\{i,\theta(i)\}$. Now $\tGamma_{\tI-\{i,\theta(i)\}}$ has two connected components interchanged by $\theta$.

\eqref{PQ} 
General building theory shows that a building is the direct product of the residues on any given chamber corresponding to the connected components of its diagram (e.g.~\cite{Ro1989a}). 
The result follows since any ${J_1}$ residue $P$ and any ${J_1^\theta}$ residue $Q^\theta$ in $R$ intersect in some chamber.

\eqref{proj}
Set $R_+=R$ and $R_-=R^\theta$.
Let $\vep\in \{+,-\}$.
First we show that $R_\vep=\proj^*_{R_\vep}(R_{-\vep})$. Namely, since $R_+$ and $R_-$ are of opposite type and contain opposite chambers, for any chamber $x_\vep\in R_\vep$ there is a chamber $x_{-\vep}\in R_{-\vep}$ opposite to $x_\vep$.
Then, the twin-apartment $\Sigma(x_+,x_-)=(\Sigma_+,\Sigma_-)$ is characterized by  $y_\vep\in \Sigma_\vep$ if and only if 
 $\delta_*(y_\vep,x_{-\vep})=\delta(y_\vep,x_\vep)$~\cite{Ti1992}.
It is coconvex~\cite{AbRo1998} and so it contains 
 $z_\vep=\proj^*_{R_\vep}(x_{-\vep})$, which is characterized by the fact that 
 \begin{align}\label{eqn:coproj}
 \delta_*(z_\vep,x_{-\vep})=\delta(z_\vep,x_\vep)=w_J,
 \end{align}
 of maximal length. Here, for any $H\sbne \tI$, $w_H$ denotes the longest word in $W_H$.
 It follows that $\delta_*(z_+,z_-)=1$ so that $\Sigma(x_+,x_-)=\Sigma(z_+,z_-)$. Hence $x_\vep=\proj_{R_\vep}^*(z_{-\vep})$ as well.
From Lemma~\ref{lem:inverse projections} we get 
 $\proj^*_{R_+}\colon R^\theta\to R$ is a (possibly type changing) isomorphism with inverse $\proj^*_{R_-}$.
 
To see how $\proj^*_{R_+}$ changes types, note that if $x'_+\in \Sigma$ is $j$-adjacent to $x_+$, for some $j\in J$ then $x'_-=\opp_\Sigma(x'_+)$ is also $j$-adjacent to $x_-$ and 
 $z'_\vep=\proj^*_{R_\vep}(x'_{-\vep})$ is $\opp_J(j)$-adjacent to $z_\vep$. Now $\opp_J$ is given by 
  \begin{align*}
  r_{\opp_J(j)}=w_J r_j w_J^{-1}.
  \end{align*} 

We have $w_J=w_{J_1}w_{{J_1^\theta}}$ and since $W_{J_1}$ and $W_{{J_1^\theta}}$ commute, we have 
\begin{align}\label{eqn:oppJ=oppKxoppK^theta}
\opp_J(j)= \begin{cases}
\opp_{J_1}(j) & \mbox{ if }j\in {J_1}\\
\opp_{{J_1^\theta}}(j) & \mbox{ if } j\in {J_1^\theta} 
 \end{cases}.
\end{align}
Thus, $\proj^*_{R_+}$ induces an isomorphism between the ${J_1^\theta}$-residue $P^\theta$
 and a ${J_1^\theta}$-residue in $R$. By~\eqref{PQ}, we may choose this residue to be $Q^\theta$.

\eqref{Q=P^theta_R}
This follows since by~\eqref{proj} $\theta_R=\proj^*_R\after \theta\colon P\to Q^\theta$ is a (possibly type-changing) isomorphism.

\eqref{AthetaR} 
Let $x=(p,q)$ with $p\in P$ and $q\in Q^\theta$.
Now $(x,x^\theta)\in R\times R^\theta$ belongs to $\De^\theta$ if and only if $(p,q)=x\opp_\De x^\theta =(p^\theta,q^\theta)$.
By~\eqref{eqn:coproj} and Lemma~\ref{lem:inverse projections}, this happens if and only if $x\opp_R x^{\theta_R}$. 
Using that $W_{J_1}$ and $W_{{J_1^\theta}}$ commute again we see that 
\begin{align*}
(p,q)\opp_R (q^{\theta_R},p^{\theta_R}) \mbox{ iff } p\opp_P q^{\theta_R} \mbox{ and }q\opp_{Q^\theta} p^{\theta_R}.
\end{align*}
By applying the isomorphism $\theta_R$, which interchanges $P$ and $Q^\theta$, we see that the latter condition is superfluous.
\epf

\ble\label{lem:tauR on Am}
Let $R$ be a residue of type $\tGamma_J\cong A_m$ for some $m$ and assume that $\proj_{R^\theta}^*$ defines a bijection between $R$ and $R^\theta$.
Then, $\theta_R$ is a type preserving automorphism of $R$.
\ele
\bpf
Note first that both $\theta$ and $\proj_{R^\theta}^*$ define a bijection between the type set of $R$ and the type set of $\theta(R)$.
Both maps can either be equal or differ by opposition.
We now prove that they cannot differ by opposition.

Let $x\in A_\theta(R)$ and consider an arbitrary twin-apartment $\Sigma$ on $x$ and $x^\theta$.
Note that $\proj_{R^\theta}^*(x)\in \Sigma$ and $\proj_R^*(x^\theta)\in \Sigma$.
Moreover, since $x\in A_\theta(R)$, the chambers $\proj_{R^\theta}^*(x)$ and $x^\theta$ are opposite in
 $R^\theta\cap \Sigma$.

Let $y=\proj^*_\pi(x^\theta)$, where $\pi$ is the $j$-panel on $x$ in $R$.
Then $y\in \Sigma\cap R$ and $l(\delta_*(y,y^\theta))=l(\delta_*(x,x^\theta))+2$ by Lemma~\ref{lem:equal distances}. More precisely, that lemma says that $y^\theta=\proj^*_{\pi^\theta}(y)$. In particular $y^\theta\in \Sigma$.

In the notation of Lemma~\ref{lem:inverse projections} $R=S$ and  so  
\begin{align*}
&l(\delta_*(x,\proj_{R^\theta}^*( x)))=l(\delta_*(y,\proj_{R^\theta}^* y)), \mbox{ and } l(\delta_*(x,x^\theta))\ne l(\delta_*(y,y^\theta)).
\end{align*} 
Therefore, by definition of projection $\delta_{-\vep}(\proj_{R^\theta}^*(y),y^\theta)\ne \delta_{-\vep}(\proj_{R^\theta}^* (x),x^\theta)=w_{\theta(J)}$.
Therefore if $\proj_{R^\theta}^* (y)$ and $\proj_{R^\theta}^*(x)$ are $j'$ adjacent, then 
 $j'$ and $\theta(j)$ are not opposite. 
\epf

\bpr\label{prop:what S is not}
Assume the terminology of Proposition~\ref{prop:AtauR}.
Then, we have the following.
\begin{enumerate}
\item  $\theta_{S_\vep}$ cannot preserve a panel,
\item $S_\vep$ cannot be of type $A_1$,
\item $S_\vep$ cannot be of type $A_2$,
\item if $S_\vep$ has type $A_1\times A_1$, then either $A_\theta(S_\vep)=S_\vep$ or $\theta_{S_\vep}$ interchanges the types.
\end{enumerate}
\epr
\bpf
Suppose $\pi$ is an $i$-panel that is preserved by $\theta_{S_\vep}$. 
Thus the bijection $\proj_{S_\vep}^*\colon S_{\vep}^\theta\to S_\vep$ restricts to a bijection between $\pi^\theta$ and $\pi$.
Note that this bijection is $\proj_{\pi}^*$.

However,  by Lemma~\ref{lem:equal distances} we see that there is a chamber $c_\vep\in \pi$ and a $w\in \delta^\theta(W)$ with the property that $\delta_*(c_\vep,c_\vep^\theta)=s_iws_{\theta(i)}$ and   
$\delta_*(d_\vep,d_\vep^\theta)=w$, for all $d_\vep\in \pi-\{c_\vep\}$ and $l(s_iws_{\theta(i)})=l(w)+2$. From the twin-building axioms it now follows that $c_\vep=\proj_{\pi}^*(d_\vep^\theta)$
 for all $d_\vep\in \pi$.
Thus, $\proj_\pi^*$ is not bijective on $\pi^\theta$, hence neither is $\proj^*_{S_\vep}$ on $S_\vep^\theta$, a contradiction.
 
Part (b) follows immediately from (a).
To see (c) note that in this case $S_\vep$ is a projective plane and any automorphism of order $2$ necessarily has a fixed point or line, hence a panel, contradicting (a).

(d) Suppose $S_\vep$ has type $A_1\times A_1$. Then,  by (a) $\theta_{S_\vep}$ cannot preserve a panel. Therefore if it fixes type, then, $\theta_{S_\vep}$ has no fixed points so that 
 $A_\theta(S_\vep)=S_\vep$.
\epf

\ble\label{lem:S a chamber}
Assume the terminology of Proposition~\ref{prop:AtauR} and set $R=R_\vep$ and $S=S_\vep$ for some $\vep=\pm$.
Suppose that $R\ne S$ and $S=A_\theta(S)$.
If $R$ has rank $2$, then $A_\theta(R)$ is connected  and if $R$ has rank $3$, then $A_\theta(R)$ is connected and simply connected.
\ele
\bpf
By Proposition~\ref{prop:AtauR},  $A_\theta(R)$ is the geometry opposite $S$.
Connectedness is proved in~\cite[Theorem 2.1]{BlBr1998},~\cite[Theorem 3.12]{Bl1999}~\cite[Proposition 7]{Ab1996}.
Now let $R$ have rank $3$.
If the diagram of $R$ is disconnected, $A_\theta(R)$ is the product of connected residues of rank $\le 2$, hence it is simply connected.
Finally suppose $R$ has type $A_3$.
If $S$ is a chamber then we are done by~\cite{Ab1996}.
In view of Proposition~\ref{prop:what S is not} this leaves the case where $S$ has type $A_1\times A_1$.
Now $A_\theta(R)$ is the geometry of all points, lines and planes of a projective $3$-space that are opposite a fixed line $l$.
That is the points and planes are those not incident to $l$ and the lines are those not intersecting $l$.
Consider any closed gallery $\gamma$ in $A_\theta(R)$.
It corresponds to a path of points and lines that all belong to $A_\theta(R)$. One easily verifies the following:
Any two points are on some plane. Hence the collinearity graph $\Xi$ on the point set of $A_\theta(R)$ has diameter $2$. Any triangle in $\Xi$ lies on a plane. 
Given any line $m$ and two points $p_1$ and $p_2$ off that line, there is a point $q$ on $m$ that is collinear to $p_1$ and $p_2$ since lines have at least three points. It follows that quadrangles and pentagons in $\Xi$ can be decomposed into triangles.
Since triangles are geometric, that is, there is some object incident to all points and lines of that triangle, $\gamma$ is null-homotopic. 
\epf

\bpr\label{prop:AtauR rank 2}
If $R$ has rank $2$, then $A_\theta(R)$ is connected.
\epr
\bpf
There are two cases: $R$ has type $A_2$ or $A_1\times A_1$.
If $R$ has type $A_2$, then by Proposition~\ref{prop:what S is not}, $S$ is a chamber and so by Lemma~\ref{lem:S a chamber} we are done.
Now let $R$ have type $A_1\times A_1$, then  $S$ is a chamber, in which case we are done again, or it is $R$.
By Proposition~\ref{prop:what S is not}, either $A_\theta(R)=R$, which is connected, or $\theta_R$ switches types and $A_\theta(R)$ is a complete bipartite graph with a perfect matching removed.
This is connected since panels have at least three elements.
\epf
\ble\label{lem:tauR disconnected}
Assume the notation of Proposition~\ref{prop:AtauR}.
Suppose that $R\cong R_1\times R_2$ and $S\cong S_1\times S_2$, where $\typ(S_i)\sbe \typ(R_i)$ for $i=1,2$.
Suppose moreover, that $\theta_S$ preserves the type sets $I_i$ of the residue $S_i$ (not necessarily point-wise).
Then, 
\begin{enumerate}
\item $\theta_R=\theta_{R_1}\times \theta_{R_2}$,
\item $A_\theta(R)\cong A_\theta(R_1)\times A_\theta(R_2)$.
\end{enumerate}
\ele
\bpf
For $i=1,2$, let $J_i=\typ(R_i)$ and let $I_i=\typ(S_i)$.
(a) 
Note that if, for $i=1,2$, $R_i'$ is a residue of type $J_i$ in $R$ then $R_1'\cap R_2'=\{c\}$ for some chamber $c$ and, for any $x\in R_1'$, $\proj_{R_2'}(x)=c$.
By assumption on $S$ the same is true for residues $S_i'$ of type $I_i$.
Note further that the same applies to the residues $R^\theta$ and 
 $S^\theta$.
Recall now that the isomorphism $R\cong R_1\times R_2$ is given by $x\mapsto (x_1,x_2)$, where $x_i=\proj_{R_i}(x)$ (see e.g.~\cite[Ch. 3]{Ro1989a}).
Thus in order to prove (a) it suffices to show that 
 \begin{align}\label{eqn:R and R_i}
 \proj_{R_i}\after \theta_R
 =\theta_{R_i}\after\proj_{R_i}.
 \end{align}
However, note that in fact
\begin{align*}
 \theta_R&=\proj_R^*\after\theta=\proj_S^*\after\theta.
\end{align*}
By Lemma 7.3 of~\cite{BlHo2009} we have
 $\proj_{S}^*=\proj_{S}^*\after\proj_{S^\theta}$ so that
 \begin{align*}
 \theta_R&=\proj_S^*\after\theta=\proj_S^*\after\proj_{S^\theta}\after\theta.
\end{align*}
The same holds for $R_i$ and $S_i$, since from~\eqref{eqn:codistance} we get 
$\proj^*_{R_i}=\proj_{R_i}\after\proj^*_R$ and 
$\proj^*_{S_i}=\proj_{S_i}\after\proj^*_S$.
Since $\theta$ is an isomorphism we also have
$\proj_{S^\theta}\after\theta=\theta\after\proj_S$, so that 
\begin{align}\label{eqn:theta R and R_i}
 \theta_R&=\proj_S^*\after\proj_{S^\theta}\after\theta=\proj_S^*\after \theta\after\proj_{S},\\
 \theta_{R_i}&=\proj_{S_i}^*\after\proj_{{S_i}^\theta}\after \theta=\proj_{S_i}^*\after\theta \after\proj_{S_i},\mbox{ for }i=1,2.\nonumber
 \end{align}
Substite~\eqref{eqn:theta R and R_i} into~\eqref{eqn:R and R_i}. For $x\in R$, 
$ \proj_{S_{i}}\after \proj_S(x)  = \proj_{S_{i}}\after\proj_{R_{i}}(x)$, and $\proj_{R_i}\after\proj_S^*=\proj_{S_i}\after\proj_S^*$, so we see that, in order to prove (a) it suffices to show that
\begin{align*}
\proj_{S_i}\after\proj_S^*\after \theta\after\proj_{S}=\proj_{S_i}^*\after \theta\after\proj_{S_i}\after\proj_S,\mbox{ for }i=1,2.
 \end{align*}
This is equivalent to showing that on $S$ we have 
 $$\proj_{S_i}\after \proj_{S}^*\after \theta=\proj_{S_i}^*\after \theta\after\proj_{S_i}, \mbox{ for }i=1,2.$$
To see this, first pick some $x\in S$ and note that if $x$ lies on the $I_2$-residue  $S_2'$,  then 
 $x,\proj_{S_1}(x)\in S_2'$, thus 
  $\theta(x),\theta\after\proj_{S_1}(x)\in S_2'^\theta$.
 But since $\theta_S$ is type-preserving, we have 
$\proj_S^*\after \theta(x),\proj_S^*\after \theta\after\proj_{S_1}(x)\in \proj_{S}^*(S_2')=S_2''$, and $S_2''$ is again of type $I_2$.
Therefore, the projection on $S_1$ of these two chambers is the same, namely $S_1\cap S_2''$.
That is, 
\begin{align*}
\proj_{S_1}\after\proj_S^*\after \theta(x)=\proj_{S_1}\after\proj_S^*\after \theta\after\proj_{S_1}(x)=S_1\cap S_2''.
\end{align*}
It is a basic property of the coprojection that $\proj_{S_1}\after\proj_S^*(y)=\proj_{S_1}^*(y)$ for any $y\in S^\theta$. Thus, we have 
\begin{align*}
\proj_{S_1}\after\proj_S^*\after\theta(x)&=
(\proj_{S_1}\after\proj_S^*)\after\theta\after\proj_{S_1}(x)\\
&=\proj_{S_1}^*\after\theta\after\proj_{S_1}(x),
\end{align*} that is, 
$\proj_{S_1}\after\theta_S=\theta_{S_1}\after\proj_{S_1}$, which proves the claim.

(b) Let $x=(x_1,x_2)\in R_1\times R_2$, and suppose $R\sbe \De_\vep$.
Then, by (a),  
 \begin{align*}
\delta_\vep(x,x^\theta)
 =&\delta((x_1,x_2),\theta_{R}(x_1,x_2))\\
 =& \delta((x_1,x_2),(\theta_{R_1}(x_1),\theta_{R_2}(x_2)))\\
 =&\delta_1(x_1,\theta_{R_1}(x_1))\cdot \delta_2(x_2,\theta_{R_2}(x_2)).
\end{align*} 
Since $A_\theta(R_1)\times A_\theta(R_2)\sbe R_1\times R_2$, 
 we see that $\delta(x,\theta_R(x))$ is maximal if and only if 
 $\delta(x_i,\theta_{R_i}(x_i))$ is maximal for $i=1,2$.
Thus $A_\theta(R)\cong A_\theta(R_1)\times A_\theta(R_2)$. 
\epf

\ble\label{lem:S=R or S A_1xA_1}
If $R$ has rank $3$, then $A_\theta(R)$ is connected and simply $2$-connected, except possibly if one of the following holds:
\begin{enumerate}
\item $R=S$, or 
\item $S<R$,  $S$ has type $A_1\times A_1$ and $\theta_S$ switches types.
\end{enumerate}
\ele
\bpf
The residue $R$ has one of three possible types: $A_3$, $A_2\times A_1$, or $A_1\times A_1\times A_1$.
By Lemma~\ref{lem:S a chamber} either $S=R$ or $S$ is a proper residue of $R$ satisfying $S\ne A_\theta(S)$.
Suppose the latter.
If $S$ is a chamber, then $S=A_\theta(S)$, which is impossible.
Moreover, by Proposition~\ref{prop:what S is not} (b) and (c), $S$ is also not a panel, or a residue of type $A_2$.
This means that $S$ has type $A_1\times A_1$ and so by Proposition~\ref{prop:what S is not} part (d), since $S\ne A_\theta(S)$, $\theta_S$ switches types on $S$.
Thus, either $S=R$, or $S$ has type $A_1\times A_1$ and $\theta_S$ switches types.
\epf

\ble\label{lem:AtauR disconnected}
Let $|\fk|\ge 3$.
If $R$ has disconnected diagram of rank $3$, then $A_\theta(R)$ is connected and simply connected.
\ele
\bpf
First suppose that $R=S$.
Then, by Corollary~\ref{cor:proj after tau has order 2}, $\theta_R=\theta_S$ has order $2$.
Whether $R$ has type $A_2\times A_1$ or $A_1\times A_1\times A_1$, the type set of $R$ can be partitioned into two non-empty sets of $\theta_R$ orbits; call them $J_1$ and $J_2$, so that $R\cong R_1\times R_2$ with $R_i$ of type $J_i$.
Taking $S_i=R_i$ for $i=1,2$, we see that  Lemma~\ref{lem:tauR disconnected} applies.
By Lemma~\ref{lem:tauR disconnected},  $A_\theta(R)\cong A_\theta(R_1)\times A_\theta(R_2)$.
By Proposition~\ref{prop:AtauR rank 2}, $A_\theta(R_i)$ is connected, hence $A_\theta(R)$ is connected and simply connected.

Next suppose that $S$ is a proper residue of $R$ of type $A_1\times A_1$ such that $\theta_S$ switches types.
As in the proof of Proposition~\ref{prop:AtauR rank 2}
 we see that $A_\theta(S)\cong S_1\times S_1^{\theta_S}-\{(x,x^{\theta_S})\mid x\in S_1\}$, for some panel $S_1$ in $S$.
 
If $R$ has type $A_1\times A_1\times A_1$, take the panel $T$ meeting $S$ in the chamber $x=S_1\cap S_1^{\theta_S}$. Then, Proposition~\ref{prop:AtauR} tells us that 
\begin{align*}
A_\theta(R)\cong&\{(t,s_1,s_2)\in T\times S_1\times S_1^{\theta_S}\mid t\not\in T\cap S, s_2\ne s_1^\theta\}\\=&(T-\{x\})\times A_\theta(S).
\end{align*} 

Since both $A_\theta(S)$ and $T-\{x\}$ are connected $A_\theta(R)$ is connected and simply connected.

We now turn to the case, where $R$ has type $A_2\times A_1$.
Let $R_i\sbe R$ be of type $A_i$ so that $R\cong R_2\times R_1$. Realize $R_2$ as the building associated to a projective plane $\Pi$ over the residue field $\fk$, representing chambers as incident point-line pairs $(p,l)$.
Identify $S_2$ with the residue in $R_2$ of a line $l_\infty$.
From Proposition~\ref{prop:AtauR} we see that $((p,l),y)\in R_2\times R_1$ belongs to $A_\theta(R)$
 iff $l\ne l_\infty$, $p\not\in l_\infty$ and $y^{\theta_S}\ne (l\cap l_\infty, l_\infty)$.
Call a point $p$ (line $l$) of $\Pi$ {\em good} if $p\not\in l_\infty$ (if $l\ne l_\infty$). Then, since $|R_1|>1$, for each chamber $(p,l)\in \Pi$ with both $p$ and $l$ good, there is a chamber $(p,l,y)\in A_\theta(R)$.
If $|\fk|\ge 3$, then to any triangle of good points and lines in $\Pi$, there is a $y\in S_2$ such that $(x,y)\in A_\theta(R)$ for any chamber $x$ on that triangle.
One verifies easily that all rank-$2$ residues meeting $A_\theta(R)$ in a chamber are connected. 
Using that all good point-line circuits $\Pi$ can be decomposed into triangles, which are all geometric, and that all rank-$2$ residues are connected we find that 
$A_\theta(R)$ is connected and simply connected.
\epf

\ble\label{lem:R=A_3}
If $R$ is of type $A_3$ and $|\fk|\ge 7$ then the geometry $A_\theta(R)$ is connected and
simply connected.
\ele
\bpf
\par{Case 1: $S=R$.}
By Lemma~\ref{cor:proj after tau has order 2}~and~\ref{lem:tauR on Am}, $\theta_R$ is an involution given by a
 semilinear map $\phi$ on a $4$-dimensional vector
space $\U$ over the residue field $\fk$. 
Since $S=R$, we also know that $\phi$ has no fixed points. Namely, the orbits of points, lines and planes have size 1 or 2; thus non-fixed points (planes) determine a fixed line and so if there is a fixed point, then either there is a fixed point-line pair or a fixed point-plane pair. However, this contradicts Proposition~\ref{prop:what S is not} (a).

Let $u,v\in \U$ be such that 
 $\phi(u)=v$ and $\phi(v)=\alpha u$ and assume $\phi$ is $\sigma$-semilinear for some $\sigma\in \Aut(\fk)$.
Then,  for any $\beta\in \fk$, we must have
\begin{align*}
\phi^2(u+\beta v)=\alpha u+\alpha^\sigma\beta^{\sigma^2} v\in \langle u+\beta v\rangle
\end{align*}
and it follows that $\alpha=\alpha^\sigma$ and $\sigma^2=1$.
Now assume that $\alpha^{-1}=\gamma\gamma^\sigma$ for some $\gamma\in \fk$, then
$\langle u+\gamma v\rangle$ is a fixed point of $\theta_R$, contradicting the previous remark.
In particular, this rules out the case where $\fk$ is finite. 

We now define the objects of the geometry $A_\theta(R)$. 
All points and all planes of $\PG(\U)$ belong to $A_\theta(R)$. The only lines in the geometry are those $2$-dimensional
spaces of $\U$ that are not
fixed by $\phi$. These will be called \textit{good lines}.
Points will be denoted by lowercase letters, good lines will be denoted by uppercase letters and planes will be denoted by greek letters.

We now describe incidence.
We shall use containment relations only for containment in $\PG(\U)$, not to be confused with incidence in $A_\theta(R)$.
Any point contained in a good line will be incident
to it and any plane containing a good line will be incident to it. A
point $p$  will be incident to a plane $\pi$ if and only if $p\sbe \pi$ and $p \not\sbe\pi^\phi\cap\pi$ (equivalently $\pi\not\spe \langle p,p^\phi\rangle$). 

We now gather some basic properties of $A_\theta(R)$.
Any two points incident to a
plane will be collinear 
and any point $p$ is incident to all planes $\pi$ so that $p\sbe \pi$
but $\pi$ does not contain the only bad line $\langle p, p^\phi\rangle$ containing $p$.
If a line $L$ is incident to a plane $\pi$, then all but one point incident to $L$ is incident to $\pi$.

Connectivity is quite immediate since any two
points $p_1, p_2$ that are not collinear will be collinear to any
other point not in the unique bad line $\langle p_1, p_2\rangle$ on $p_1$ (and $p_2$).

In order to prove simple connectivity we first reduce any path to a
path in the collinearity graph. Indeed any path $p_1\pi p_2$ will be
homotopically equivalent to the path $p_1Lp_2$ where $L=\langle
p_1,p_2\rangle$. 
Any path $p\pi L$ will be homotopically
equivalent to the path $pL'p'L$ where $p'$ is a point on $L$ that is also incident to $\pi$ and $L'=\langle p,p'\rangle$. 
Note that since $p'$ is incident to $\pi$, $L'$ is a good line.
Finally a path
$L_1\pi L_2$ is homotopically equivalent to the path $L_1p_1L'p_2L_2$
where $p_i$ are points on $L_i$ that are incident to $\pi$ and
$L'=\langle p_1,p_2\rangle$.

Therefore, to show simple connectedness we can restrict to paths in the collinearity graph. Note
also the fact that if $p$ is a point and $L$ is a good line not
incident to $p$ then $p$ will be collinear to all but at most one
point on $L$ (namely the intersection of the unique bad line on $p$ and $L$ if this intersection exists).  This enables the decomposition of any path in the collinearity graph to triangles. Indeed, the diameter of the collinearity graph is two and so any path can be decomposed into triangles, quadrangles and pentagons. Moreover, if $p_1,p_{2},p_{3},p_{4}$ is a quadrangle then, since $|k|\ge 4$, the line $\langle p_{2}, p_{3}\rangle$ will admit a point collinear to both $p_{1}$ and $p_{4}$ decomposing the quadrangle into triangles. Similarly, if  $p_1,p_{2},p_{3},p_{4}, p_{5}$ is a pentagon,  then there will be a point on the good line $\langle p_{3}, p_{4}\rangle$ that is collinear to $p_{1}$. Thus, the pentagon decomposes into quadrangles.  Therefore it suffices to decompose triangles into geometric triangles.

Assume that $p_1, p_2, p_3$ is a triangle. The plane $\pi=\langle
p_1,p_2, p_3\rangle$ is incident to all three (good) lines in the triangle and so, either the
triangle is geometric and then we are done, or one of the points is not incident to
$\pi$, that is, it lies on the bad line $\pi\cap\pi^\phi$. Since the triangle lines are good, there is at most one such point. Let us assume that $p_1$ is not incident to $\pi$.

Consider a plane $\pi'$ that contains the line $\langle
p_2,p_3\rangle$ and so that $p_2$ and $p_3$ are incident to
$\pi'$. This is certainly possible since $|\fk|\ge 4$ and one only need to stay clear
of the planes $\langle p_2,p_3, p_3^\phi\rangle$ and $\langle
p_2,p_3, p_2^\phi\rangle$.

Note that by choice of $\pi'$, any line $L$ with $p_i\sbe L\sbe \pi'$ ($i=2,3$) is good.
Let now for each $i=2,3$
\begin{align*}
\cL_i=\{L \mbox{ is a good line in }\pi' \mid p_i\sbe L, p_1, p_i \mbox{ are incident to }\langle p_1, p_i,L\rangle\}.
\end{align*}
 The only lines of $\pi'$ on $p_i$ not in $\cL_i$ 
 are $\langle p_2,p_3\rangle$ 
and $\langle p_1,p_i, p_i^\phi\rangle \cap \pi'$ so $\cL_i=|\fk|-1$.
Note that if $L\in \cL_i$ then 
the only point incident to $L$ not incident to $\pi'$ is $L\cap \pi'\cap{\pi'}^\phi$.
Pick distinct lines $L_{i,j}\in \cL_i$ with  $j=1,2,3,4$.
Of the $16$ intersection points $p_{j,k}=L_{2,j}\cap L_{3,k}$ at most $8$ are not incident to 
 one of the three planes that they define.
For instance, each of the four planes $\langle p_1, L_{2,j}\rangle$ contains exactly one bad line.
 This bad line can be on at most one of the four intersection points $p_{j,k}$ $k=1,2,3,4$.
Thus, there must be at least $16-8=8$ points $p_{j,k}$ that are not incident to the bad lines in $\langle p_1,L_{2,j}\rangle$ or $\langle p_1, L_{3,k}\rangle$.
Out of these $8$ points, at most four are on the bad line $\pi'\cap \pi'^\phi$. 
Using any of the remaining $4$ points $p$, the triangle 
$p_1p_2p_3$ can be decomposed into the geometric triangles consisting of $p$ and two points from $\{p_1,p_2,p_3\}$.

\par{\bf Case 2: $S$ of type $A_1\times A_1$ and $\theta_S$ switches types.}
 The geometry is rather
similar to the previous one. There is a line $\mathbf{L}$ so that $S$ is the
residue corresponding to $\mathbf{L}$ and the map $\theta_S$ induces a pairing
between points of $\mathbf{L}$ and planes on $\mathbf{L}$.
The geometry $A_{\theta}(R)$ is
described as follows. The points of the geometry are all the points of
$\U$  not in $\mathbf{L}$, the
lines of the geometry are all the lines of $\U$ not intersecting $\mathbf{L}$ and the planes are all
planes of $\U$ not containing $\mathbf{L}$. 

We now describe incidence.
Any line included in a plane is incident to
it and any point included in a line is incident to it. A point $p$ is
incident  to a plane $\pi$ if and only if the plane $\pi'=\langle
p,\mathbf{L}\rangle$ is not paired to the point $p'=\mathbf{L} \cap \pi$; that is $\pi'^\phi\ne p'$.

We now gather a few useful properties of this geometry.
Note a
number of similarities with the previous geometry. Any plane $\pi$ is
incident to all the points $p\sbe \pi$  that are not contained in the unique bad line $\lambda(\pi)=\pi'  \cap \pi $ on $\pi$; here $\pi'$ is the
plane paired to the point $\pi \cap \mathbf{L}$. Dually  any point $p$ is
incident to all the planes $\pi\spe p$ that do not contain the unique bad line $\lambda(p):= \langle p,
p'\rangle$ on $p$; here
$p'$ is the point paired to the plane $\langle p, \mathbf{L} \rangle$. If $p$ is a point and $L$ is a good line not
incident to $p$ then $p$ will be collinear to all but one
point on $L$; namely the non-collinear point on $L$ is the intersection of $L$ with the bad plane $\langle p,{\mathbf L}\rangle$.

 Any two
points $p_1, p_2$ that are not collinear have the property that  
 $\langle p_1,p_2\rangle $ intersects ${\mathbf L}$ and so 
  any point not in $\langle p_1,p_2,{\mathbf L}\rangle$ will be collinear to both $p_1$ and $p_2$.
In particular, the geometry $A_\theta(R)$ is connected and the diameter of the collinearity graph is $2$.

The reduction to the collinearity graph is a little more involved
because not every two points on a good plane will be collinear. However
any two non-collinear points incident to a good plane $\pi$ are collinear to any other
point $p_3$ incident to $\pi$ but not in the line $p_1p_2$ since ${\mathbf L}$ intersects $\pi$ in exactly one point.

The previous remark immediately shows that a path of type 
 $p_1\pi p_2$ can be replaced by a path $p_1, L_1,p',L_2,p_2$, where all elements are incident to $\pi$.
Suppose we have a path of type $p\pi L$.
Since $\pi$ is incident to all but one point on the line $L$ and 
 $p$ is collinear to all but one point on the line $L$, we can replace this path by one of type $pL_1p_2L$, where
  all objects are incident to $\pi$.
Suppose we have a path of type $L_1\pi L_2$.
This reduces to the previous case since all but one point of $L_1$ are incident to $\pi$.

As before, given any line $L$ and two points $p_1$ and $p_2$ not on $L$, there are only two points on $L$ that are not collinear to at least one of $p_1$ and $p_2$.
The proof that all paths in the collinearity graph decompose into
triangles is identical. Therefore it suffices to show that any triangle decomposes into geometric triangles.

We now modify the argument above to decompose
triangles. 
Again our aim is to select a point $p_0$ not on $\pi$ collinear to $p_k$ ($k=1,2,3$), and such that 
$p_0$, $p_i$, and $p_j$ are incident to $\pi_{i,j}=\langle p_0,p_i,p_j\rangle$ ($1\le i<j \le 3$).
The only difference is once more the fact that two points incident to a good plane are collinear if and only if the line joining them does not pass through $\mathbf{L}$. 

To ensure that $p_2$ and $p_3$ are incident to $\pi'=\pi_{2,3}$, let $\pi'$ be a plane on $p_2p_3$ that does not contain $\lambda(p_2)$ or $\lambda(p_3)$.
Let now for each $i=2,3$
\begin{align*}
\cL_i=\{L\mbox{ is a good line in }\pi'\mid p_i\sbe L, p_1, p_i \mbox{ are incident to }\langle p_1, L\rangle\}.
\end{align*}
Since each $L\in \cL_i$ is good, any $p_0\sbe L$ is collinear to $p_i$. Moreover $\pi_{1,i}=\langle p_1,L\rangle$.

In order to ensure that $p_1$ is incident to $\langle p_1, L\rangle$,  $\langle p_1, L\rangle$ must not contain $\lambda(p_1)$, that is we must exclude $p_2p_3$ from $\cL_i$.
In order to ensure that $p_i$ is incident to $\langle p_1, L\rangle$ we must exclude the line $\langle p_1,\lambda(p_i)\rangle\cap \pi'$ from $\cL_i$ ($i=2,3$).
Let $\mathbf{p'}=\mathbf{L}\cap \pi'$.
To ensure that $p_0$ and $p_i$ are collinear we must exclude the line $p_i \mathbf{p'}$ from $\cL_i$.
As a consequence the sets $\cL_i$ have $|\fk|-2$ lines.

Now assume that $|\fk|\ge 7$. Then, pick lines
 $L_{i,j}\in \cL_i$ ($i=2,3$, $j=1,2,3,4,5$) and define the set $\sP=\{L_{2,i}\cap L_{3,j} \mid i,j =1,2,3,4,5\}$ of size $25$. 
 Note that if $p_0 \in \sP$ then $p_0$ is collinear to $p_2$ and $p_3$, $p_1,p_i$ are incident to $\pi_{1,j}$. We still need to insure that $p_0$ is collinear to $p_1$ and $p_0$ is incident to $\pi_{i,j}$.
 
In order to ensure that $p_0$ is collinear to $p_1$, we must choose $p_0$ so that $p_0p_1$ does not intersect $\mathbf{L}$.This means that $p_0$ does not lie on the line $\langle p_1,\mathbf{L}\rangle\cap \pi'=\langle \mathbf{p'},(\lambda(p_1)\cap p_2p_3)\rangle$. This eliminates at most the 5 points $L_{2,j}\cap \langle p_1,\mathbf{L}\rangle\cap \pi'$ from $\sP$.

 To ensure that $p_0$ is incident to $\pi_{2,3}=\pi'$ we must choose $p_0$ off $\lambda':=\lambda(\pi')\spe \mathbf{p'}$. This eliminates at most the 5 points $L_{2,j}\cap \lambda'$  from $\sP$.
 
Finally in order to ensure that $p_0$ is incident to $\pi_{1,i}$, we notice that each of the 10 planes  $\pi_{1,i}=\langle p_1,L_{i,j}\rangle$ has a unique bad line and so at most one of the points of $L_{i,j}$ fails to be incident to this plane. This eliminates at most 10 more points  from $\sP$.
If $p_0\in \sP$ is any of the remaining points, of which there are at least $5$,  then  $p_0,p_i,p_j$ are all geometric triangles. This decomposes the initial triangle  $p_1,p_2,p_3$ into geometric triangles.
\epf
\bth\label{thm:AtauR rank 3}
Suppose that $|\fk|\ge 7$.
If $R$ has rank $3$, then $A_\theta(R)$ is connected and simply $2$-connected.
\eth
\bpf
The theorem follows from Lemmas~\ref{lem:S=R or S A_1xA_1},~\ref{lem:AtauR disconnected}~and~\ref{lem:R=A_3}.
\epf

\subsection{Proof of Theorem~\ref{mainthm:omega kac-moody construction}}\label{subsec:proof of theorem 3}
In order to prove Theorem~\ref{mainthm:omega kac-moody construction}, we first note that $\amG^\delta\cong\amL^\delta$.
This follows from Proposition~\ref{prop:non-oriented CT amalgam}.

We shall now prove the theorem using Proposition~\ref{prop:efficient completion recognition}. 

For $\tJ\sbne \tI$ 
 and $\vep=+,-$, let $R_{\tJ,\vep}$ be the $\tJ$-residue of $\De_\vep$ on $c_\vep$.  
Also, let $K_{\tJ}$ be the Levi component of the standard parabolic subgroup in $\G$ stabilizing the pair $(R_{\tJ,+},R_{\tJ,-})$.
Now let $J\sbne I$ and by abuse of notation view $I\sbe \tI$, and let $\tJ=J\cup J^\theta$.
Write $R_\tJ=R_{\tJ,+}$, then, by Proposition~\ref{prop:tau on Delta}, $R_{\tJ,-}=R_{\tJ}^\theta$, and $(R_\tJ,R_\tJ^\theta)$ intersects $(\De^\theta,\approx)$ in a  residue of $\De^\theta$.

Let  $B_J$ be the stabilizer in $\G^\theta$ of the residue $(R_\tJ,R_\tJ^\theta)\cap \De^\theta$.
Then, 
\begin{align*}
\amB=\{B_J\mid J\sbne I\}
\end{align*}
with connecting maps given by inclusion of subgroups in $\G^\theta$, is the amalgam of maximal parabolic subgroups of $\G^\theta$ for the action on $\De^\theta$. 
Recall from Proposition~\ref{prop:tau on Delta} that $\theta(m)=m-n \mod{2n}$, for $m\in I=\{1,2,\ldots,n\}$.
For $m\in \{1,\ldots,n\}$ write $B_{(m)}=B_{I-\{m\}}$.

\ble\label{lem:G theta completion of B}
The universal completion of the amalgam
 $\amB$ equals $\G^\theta$.
\ele
\bpf
Under the assumptions of Theorem~\ref{mainthm:omega kac-moody construction}, $n\ge 4$, and $|\fk|\ge 7$, so that 
by Proposition~\ref{prop:AtauR rank 2} and Theorem~\ref{thm:AtauR rank 3} the residual filtration $\cC$ satisfies the conditions of Theorem~\ref{thm:DeMu}, noting that by Proposition~\ref{prop:A_tau has unique distance}, $\aff(R)=A_\theta(R)$.
It follows that $(\cC_0,\sim)\cong (\De^\theta,\sim)$ is connected and simply connected and hence by Lemma~\ref{lem:sim to approx}, so is $(\De^\theta,\approx)$. 
As mentioned above, since $\fk/\fk_\kaut$ is cyclic and Galois, Theorem~\ref{thm:flag-transitive} tells us that $\G^\theta$ is a flag-transitive automorphism group of $\De^\theta$.
Therefore, by Tits' Lemma~\cite[Corollaire 1]{Ti1986b}, $\G^\theta$ is the universal completion of the amalgam $\amB$. 
\epf

\mn
Recall that $\amL^\delta=\{\L_{i},\L_{ij}\mid i,j\in \{1,2,\ldots,n\}\}$, with $\L_{i}$ and $\L_{ij}$ as defined in~Subsection~\ref{subsec:amL^delta}.
For $\emptyset\sbne J\sbne I$, let
 \begin{align*}
 \L_J=\langle \L_{i},\L_{i,j}\mid i,j\in J\rangle_{\G^\theta}.
 \end{align*}
Recall from Definition~\ref{dfn:standard pair} that, for each $m\in \{1,2,\ldots,n\}$,  $D_m$ denotes the diagonal torus of $\G_m\in \amG^\delta$.
As $\amG^\delta$ has property (D), we may identify  $D_m$ unambiguously with its image $\L_m\cap \Dt$ in $\G^\theta$. Let $K_\tJ^\theta=K_\tJ\cap \G^\theta$.
\bpr\label{prop:theta parabolics}
In the notation from this subsection, we have 
\begin{enumerate}
\item\label{parabolics} 
$B_J=K_\tJ^\theta$,
\item\label{decomp}
$B_J=\langle \L_J,\Dt\rangle_{\G^\theta}$.
\end{enumerate}
\epr
\bpf
\eqref{parabolics}
Clearly, $K_\tJ\cap \G^\theta\le B_J$. 
Conversely, for $g\in B_J$, in view of~\eqref{eqn:A_theta(R) in De^theta}, we have $A_\theta(R_\tJ)\sbe g(R_\tJ)\cap R_\tJ$, but since $R_\tJ$ and $g(R_\tJ)$ have the same type, they must be equal and the same holds for $R_\tJ^\theta$.
Hence, in fact $B_J= K_\tJ\cap \G^\theta$.

\eqref{decomp}
Let $J\sbe I$ and let $J=\cup_{i}J_{i}$ be a decomposition of $J$ corresponding to connected components of the diagram $\Gm_J$ induced on the node set $J$.
If necessary using $\Phi_{{\R_{\kaut^2}},2n}$, we may assume that $n\not\in J$.

Now $\L_J\le \G^\theta$ stabilizes the $\tJ$-residue of $\De_\vep$ on $c_\vep$, so that
 $\L_J\le K_\tJ^\theta$.
Also, by Lemma~\ref{lem:GUD},  $\Dt= \SUD\le K_\tJ^\theta$.
Thus, by \eqref{parabolics} $B_J\ge \L_J\Dt$.

To see the reverse inclusion, first note that 
since $\G^\theta$ is flag-transitive on $\De^\theta$, $B_J$ is transitive on the chambers of  $(R_\tJ,R_\tJ^\theta)\cap \De^\theta$. By definition, $B_J$ contains $B_\emptyset=\SUD=\Dt$, which is the stabilizer of $(c_+,c_-)$ in $\G^\theta$.

We now show that $\langle \L_J,\Dt\rangle$ has these same properties.
Since this group is a subgroup of $B_J$ it acts on $A_\theta(R_\tJ)$. By Proposition~\ref{prop:theta residues}, $A_\theta(R_\tJ)$ consists of the pairs of chambers $(p,q^{\theta_R})$ in $P\times P^{\theta_R}$,  where $(p,q)$ is a pair of opposite chambers in the building $P$ whose diagram is the subdiagram of $\tA_{2n-1}$ induced on $J$.
Then, $\L_{J}$ acts as $\oplus_{i} \SL_{n_{i}+1}(k)$ (where $|J_{i}|=n_i$) on $P$. 
Therefore, it is certainly transitive on pairs of opposite chambers in $P$. Thus, $\langle \L_J,\Dt\rangle$ is transitive on the chambers of $(R_\tJ,R_\tJ^\theta)\cap \De^\theta$.
Moreover, this group contains $\Dt$, which coincides with $\SUD$ by Lemma~\ref{lem:GUD}. We are done.
\epf

\mn
In the notation of Proposition~\ref{prop:efficient completion recognition}, Proposition~\ref{prop:theta parabolics} demonstrates that the amalgam $\amB$ is indeed the amalgam $\amB$ as constructed in~\eqref{eqn:amB}, Lemma~\ref{lem:G theta completion of B} proves that $\amB$ and $\G^\delta$ satisfy condition (c), 
 and Proposition~\ref{prop:non-oriented CT amalgam} shows that $\amL^\delta$ satisfies condition (a).

Therefore it remains to show that condition (b) of Proposition~\ref{prop:efficient completion recognition}
  is satisfied.
\ble
The group $H(\G^\delta)$, as defined in Lemma~\ref{lem:H} is trivial.
\ele
\bpf
This follows by noting that if $a=a^\delta=a^{-\kaut}$ ($\delta=\kaut\tau$), then taking the product over all $\phi_i$ images of the matrix
 \begin{align*}
 d(a)=\begin{pmatrix}
 a & 0\\
 0 & a^{-1}
 \end{pmatrix}
 \end{align*}
we obtain the identity of $\SL_{2n}(\fk)$.
Indeed 
\begin{align*}
 \prod_{i=1}^{n-1}\phi_i(d(a))
&=\left(\begin{array}{@{}ccc|ccc@{}}
a &  &  & &&\\
 & I_{n-2} &&& & \\
 &  &  a^{-1} & & & \\
 \hline
 &  & & a^{-\kaut} & &  \\
 &  & &  & I_{n-2}  &  \\
  & & & && a^{\kaut}\\
  \end{array}
  \right)
  \mbox{ and }\\
    \phi_n(d(a))
&=\left(\begin{array}{@{}ccc|ccc@{}}
a^{\kaut^{-1}} &  &  & &&\\
 & I_{n-2} &&& & \\
 &  &  a & & & \\
 \hline
 &  & & a^{-1} & &  \\
 &  & &  & I_{n-2}  &  \\
  & & & && a^{-\kaut}\\
  \end{array}
  \right).
\end{align*}
\epf

\mn
Theorem~\ref{mainthm:omega kac-moody construction}
 now follows from Proposition~\ref{prop:efficient completion recognition}.
 
The conditions $n\ge 4$ and $|\fk|\ge 4$ come from the classification result in~\cite{BloHof2013}. The condition $|\fk|\ge 7$ is used to show connectedness and simple connectedness of $\De^\theta$ (Theorem~\ref{thm:AtauR rank 3}), the condition that $\fk/\fk_\kaut$ be cyclic and Galois ensures that $\G^\theta$ is flag-transitive on $\De^\theta$ (Theorem~\ref{thm:flag-transitive}), and the condition that $N_{\fk_{\kaut^2}/\fk_\kaut}$ is surjective is used to show that $\Dt$ is the full stabilizer in $\G^\theta$ of a pair of opposite chambers in $\De^\theta$ (Lemma~\ref{lem:GUD}).

\end{document}